\documentclass[smallcondensed,english]{svjour3}
\usepackage[T1]{fontenc}
\usepackage[11pt]{extsizes} 
\usepackage[latin9]{inputenc}
\usepackage{color}
\usepackage{babel}
\usepackage{verbatim}
\usepackage{float}
\usepackage{url}
\usepackage{amssymb}
\usepackage{amstext}
\usepackage{graphicx}
\usepackage[numbers]{natbib}
\usepackage[unicode=true,pdfusetitle,
 bookmarks=true,bookmarksnumbered=false,bookmarksopen=false,
 breaklinks=false,pdfborder={0 0 1},backref=false,colorlinks=false]
 {hyperref}
\usepackage{breakurl}

\makeatletter

\providecommand{\tabularnewline}{\\}
\floatstyle{ruled}
\providecommand{\algorithmname}{Algorithm}
\floatname{algorithm}{\protect\algorithmname}

 \usepackage{fullpage, graphicx,psfrag,amsmath,amsfonts,verbatim}
\usepackage[small,bf]{caption}
\usepackage{algorithm}
\usepackage{algpseudocode}
\usepackage{mathtools}
\newtheorem{assumption}{Assumption}
\usepackage{ dsfont }
\date{}

\makeatother

\begin{document}

\title{On seeking efficient Pareto optimal points in multi-player minimum
cost flow problems with application to transportation systems\thanks{A preliminary version of this work without proofs appeared  in \cite{gupta2016multi}.}}

\author{Shuvomoy Das Gupta \& Lacra Pavel}

\institute{ Shuvomoy Das Gupta \at                
Research \& Technology Department\\
Thales Canada, Transportation Solutions \\  
105 Moatfield Drive, Toronto, Ontario, Canada \\              
\email{shuvomoy.dasgupta@thalesgroup.com} 
\and   Lacra Pavel \at       
 Systems Control Group  \\         
Department of Electrical \& Computer Engineering\\
 University of Toronto  \\
10 King's College Road, Toronto, Ontario, Canada \\               
\email{pavel@control.utoronto.ca}
   }
\maketitle
\begin{abstract}
\textcolor{black}{In this paper, we propose a multi-player extension
of the minimum cost flow problem inspired by a transportation problem
that arises in modern transportation industry. We associate one player
with each arc of a directed network, each trying to minimize its cost
function subject to the network flow constraints. In our model, the
cost function can be any general nonlinear function, and the flow
through each arc is an integer. We present algorithms to compute }\textcolor{black}{\emph{efficient
Pareto optimal point(s}}\textcolor{black}{), where the maximum possible
number of players (but not all) minimize their cost functions simultaneously.
The computed Pareto optimal points} are Nash equilibriums\textcolor{black}{{}
if the problem is transformed into a finite static game in normal
form}. 
\end{abstract}

\section{Introduction}

\paragraph*{}

In recent years, product transportation systems
are increasingly being dominated by retailers such as Amazon, Alibaba, and
Walmart, who utilize
e-commerce solutions to fulfill customer supply chain expectations
\citep{turban2017electronic,galloway,li2014managing,mendez2017beyond}. In the supply chain strategy of these retailers,
products located at different warehouses are shipped to geographically
dispersed retail centers by different transportation organizations.
These transportation organizations (\emph{carriers}) compete among themselves and transport goods between warehouses and
retail centers over multiple transportation links. For example, Amazon
uses FedEx, UPS (United Parcel Service), AAR (Association of American
Railroads), and other competing organizations to provide transportation
services \citep[Chapter 11]{laseter2011internet}. 

Product shipment from a warehouse to a retail center
requires contracting multiple competing carriers, \emph{e.g.}, a common shipment
may comprise of (i) a trucking shipment from the warehouse to a railhead
provided by FedEx, then (ii) a rail shipment provided by AAR, and
finally (iii) a trucking shipment from the rail yard to the retail
center provided by UPS. It is common that different competing carriers
operate over the same transportation link, \emph{e.g.}, both FedEx and UPS provide
trucking shipment services for Amazon over the same transportation
link \citep[Section 9.1]{li2014managing}. The goal of each carrier
is to maximize its profit (minimize its cost). So a relevant question
in this regard is how to determine a good socially-optimal solution
for the competing carriers. 

We can formulate the multi-carrier transportation
setup described above as a multi-player extension of the well-known
minimum cost flow problems. The transportation setup can be modeled
by a directed network, where a warehouse is a supply node and a retail
center is a demand node. A transportation link is an arc, and a carrier (\emph{e.g.}, FedEx, UPS)
in charge of transporting products over that transportation link is
a player. The products transported over the directed network are modeled
as the flow through the network, and customer supply chain expectations
can be modeled as mass-balance constraints. Each of the carriers is
trying to maximize its profit by maximizing the total number of products
that it transports. Carriers are competing for a limited resource,
namely the total number of products to be transported. Note that one
carrier making the maximum profit may impact the other carriers in
a negative manner and even violate customer supply chain expectations.
Our goal is to define a socially-optimal solution concept in this
setup and to provide algorithms to calculate such a solution.

The problem above can be generalized as a multi-player
minimum cost flow problem \citep{gupta2016multi}. We associate one
player with one arc of the directed and connected network graph. Parallel
arcs between two nodes denote two competing carriers over the same
transportation link. Each of the players is trying to minimize its
cost function, subject to the network flow constraints. The flow through
each arc is taken to be an integer. This assumption incurs no significant
loss of generality because by suitable scaling we can use the integer
model to obtain a real-valued flow vector to any desired degree of
accuracy. Naturally, defining an efficient solution concept is of
interest in such a problem setup.

\textcolor{black}{The unambiguously best choice would be a utopian
}\textcolor{black}{\emph{vector optimal solution,}}\textcolor{black}{{}
which minimizes all the objectives simultaneously. However, this is
unlikely to exist in practice \citep[page 176]{Boyd2004}. A}\textcolor{black}{\emph{
generic Pareto optimal point, }}\textcolor{black}{where none of the
objective functions can be improved without worsening some of the
other objective values, is a better solution concept. However, there
can be numerous such generic Pareto optimal points, many of which
would be poor in quality or efficiency \citep[Section 4.7.5]{Boyd2004}.
In this paper, we investigate an }\textcolor{black}{\emph{efficient
Pareto optimal point}}\textcolor{black}{{} \citep[Section 2.2]{Miettinen2012}
as a good compromise solution that finds a balance between the utopian
vector optimality and the generic Pareto optimality. We present algorithms
to compute an efficient Pareto optimal point, }\textcolor{black}{\emph{i.e.},}\textcolor{black}{{}
a Pareto optimal point where the maximum possible (but not all) number
of players minimize their cost functions simultaneously.}

\paragraph*{Related work.} $\;$ The classic version of the minimum cost flow problem has a linear cost function
for which polynomial time algorithms exist \citep[Chapter 10]{Ahuja1988},
even when the network structure (\emph{e.g.}, nodes and arcs) is subject
to uncertainty \citep{Bertsimas2013}. The polynomial runtime can
be improved to a linear runtime when the minimum cost flow problem has a special
structure \citep{Vaidyanathan2010}. However, for nonlinear cost functions,
results exist for very specific cases, and no work seems to exist
for arbitrary nonlinear functions. The most commonly used nonlinear
cost function is the fixed charge function, where the cost on the
arc is $0$ if the flow is zero and affine otherwise. Network flow problems
with fixed charge costs are studied in \citep{Magnanti1984,Graves1985,Daskin2011},
where the integrality condition on the flow is not considered. Minimum cost flow problems with concave cost functions are
studied in \citep{Feltenmark1997,Yaged1971,He2015,Tuy1995,Fontes2006}. Minimum cost flow problems with piece-wise linear cost functions is investigated in \citep{Zangwill1968,Jorjani1994,Gamarnik2012}. A dynamic domain contraction algorithm for nonconvex piece-8wise linear network flow problems is proposed in \cite{Kim2000}. A particle swarm optimization based hybrid algorithm to solve the minimum concave cost network flow problem is investigated in \cite{Yan2011}. Finding Pareto optimal points in multi-objective network flow problem
with integer flows has been limited so far to linear cost functions
and two objectives \citep{Raith2009,Lee1993,Eusebio2009, Hernandez2011}. A multi-player
minimum cost flow problem with nonconvex cost functions is explored in \citep{gupta2016multi}. Integer multi-commodity flow problems
are investigated in \citep{Barnhart1996,Brunetta2000,Ozdaglar2004}; the underlying problems are optimization problems in these papers.  
\begin{comment}
 and the solution concepts are optimal solutions to them. In comparison,
in our paper we have focused on a multi-player extension of such problems
for a single commodity and have focused on computing efficient Pareto
optimal points as solution concepts.
\end{comment}

\paragraph*{Contributions.} $\;$ In this paper, we propose an extension of the minimum
cost flow problem to a multi-player setup and construct algorithms
to compute efficient Pareto optimal solutions. Our problem can be
interpreted as a multi-objective optimization problem \citep[Section 4.7.5]{Boyd2004}
with the objective vector consisting of a number of univariate general
nonlinear cost functions subject to the network flow constraints and
integer flows. In comparison with existing literature, we do not require
the cost function to be of any specific structure. The only assumption
on the cost functions is that they are proper. In contrast to relevant
works in multi-objective network flow problems, our objective vector
has an arbitrary number of components; however, each component is a
function of a decoupled single variable. We extend this setup to a problem 
class that is \emph{strictly} larger than, but that contains the
network flow problems. We develop our algorithms for this larger class. 

We show that, although in its original form the problem has coupled
constraints binding every player, there exists an equivalent variable
transformation that decouples the optimization problems for a maximal
number of players. Solving these decoupled optimization problems can
potentially lead to a significant reduction in the number of candidate
points to be searched. We use the solutions of these decoupled optimization
problems to reduce the size of the set of candidate efficient Pareto
optimal solutions even further using algebraic geometry. Then we present
algorithms to compute efficient Pareto optimal points that depend
on a certain semi-algebraic set being nonempty. We also present a penalty based approach applicable when that
semi-algebraic set is empty; such an approach can be of value to network administrators and policy makers.
To the best of our knowledge, our methodology is novel. The computed
efficient Pareto optimal point has some desirable properties: 
\textcolor{black}{(i)}\textcolor{black}{{} it is a good compromise
solution between the utopian vector optimality and the generic Pareto
optimality, and (ii) it is a Nash equilibrium if we convert our setup into
a finite static game in normal form.}

The rest of the paper is organized as follows. Section \ref{sec:Notation-and-notions}
presents notation and notions used in the paper. Section \ref{sec:Prb-stmnt-for-ntworks}
describes the problem for directed networks and the extension to a
strictly larger class. In Section \ref{sec:transforming_the_game},
we transform the problem under  consideration into decoupled optimization
problems for a number of players and reformulate the optimization problems
for the rest of the players using consensus constraints. Section \ref{sec:Algorithms-to-compute}
presents algorithms for computing efficient Pareto optimal points
for our problem if a certain semi-algebraic set is nonempty. Section
\ref{sec:A-penalty-based} discusses a penalty based approach if the
semi-algebraic set is empty. Section \ref{sec:Numerical-Examples}
presents an illustrative numerical example of our methodology in a
transportation setup. Finally, in Section \ref{sec:Conclusion}
we present some concluding remarks regarding our methodology. \textbf{Proofs
are provided in the Appendix.}

\paragraph{Notation and notions. \label{sec:Notation-and-notions}}

We denote the sets of real numbers, integers, and natural numbers by
$\mathbf{R}$, $\mathbf{Z}$, and $\mathbf{N}$, respectively. The $i$th
column, $j$th row, and $(i,j)$th component of a matrix $A\in\mathbf{R}^{m\times n}$
is denoted by $A_{i}$, $a_{j}^{T}$, and $a_{ij}$, respectively. The submatrix
of a matrix $A\in\mathbf{R}^{m\times n}$, which constitutes of its
rows $r_{1},r_{1}+1,\ldots,r_{2}$ and columns $c_{1},c_{1}+1,\ldots,c_{2}$,
is denoted by $A_{[r_{1}:r_{2},c_{1}:c_{2}]}\in\mathbf{R}^{(r_{2}-r_{1}+1)\times(c_{2}-c_{1}+1)}$.
If we make two copies of a vector $x\in\mathbf{R}^{n}$, then the copies are denoted by $x^{(1)}$ and $x^{(2)}$. By $I_{i,j}\in\mathbf{R}^{n\times n}$,
we denote a matrix that has a $1$ on its $(i,j)$th position and
$0$ everywhere else. If $C$ and $D$ are two nonempty sets, then $C+D=\{x+y\mid x\in C,y\in D\}$.
If $A\in\mathbf{R}^{m\times n}$ is a matrix and $C$ is a nonempty
set containing $n$-dimensional points, then $AC=\{Ax\mid x\in C\}$.
\textcolor{black}{The set of consecutive integers from 1 to $n$ is
denoted by $[n]=\{1,2,\ldots,n\}$ and $m$ to $n$ is denoted by
$[m:n]=\{m,m+1,\ldots,n\}$.} If we have two vectors $x,y\in\mathbf{R}^{n}$,
then $x\succeq y$ means 
\[
\left(\forall i\in[n]\right)\quad x_{i}\geq y_{i},
\]
 and we write $x-y\in\mathbf{R}_{+}^{n}$. Depending on the context, $0$ can represent the scalar zero, a column vector of zeros, or a matrix with all the entries zero,\emph{ e.g.}, if we say $x \in \mathbf{R}^n$ and $x=0$, then $0$ represents a column vector of $n$ zeros. On the other hand, if we say $A \in \mathbf{R}^{m \times n} $ and $A=0$, then $0$ represents a matrix of $m$ rows and $n$ columns with all entries zero.

Consider a standard form polyhedron $\{x\in\mathbf{R}^{n}\mid Ax=b,x\succeq0\}$,
where $A\in\mathbf{R}^{m\times n}$ is a full row rank matrix, and
$b\in\mathbf{R}^{m}$. A \emph{basis matrix} of this polyhedron is
constructed as follows. We pick $m$ linearly independent columns
of $A$ and construct the $m\times m$ invertible square submatrix
out of those columns; the resultant matrix is called a basis matrix.
The concept of basis matrix is pivotal in simplex algorithm, which
is used to solve linear programming problem. Suppose we have a polyhedron
defined by linear equality and inequality constraints in $\mathbf{R}^{n}$. Then $\tilde{x}\in\mathbf{R}^{n}$ is called a \emph{basic solution} of the polyhedron, if  all
the equality constraints are active at $\tilde{x}$, and out of all the active constraints
(both equality and inequality) that are active at $\tilde{x}$, there
are $n$ of them that are linearly independent. 

\section{\label{sec:Prb-stmnt-for-ntworks}Problem statement}

\label{subsec:Problem-statement-for-networks}

\paragraph*{}

\textcolor{black}{Let $G=(\mathcal{\mathcal{M}},\mathcal{A})$ be
a }\textcolor{black}{\emph{directed connected}}\textcolor{black}{{}
graph associated with a network, where $\mathcal{M}=[m+1]$ is the
set of nodes, and $\mathcal{A}=[n]$ is the set of (directed) arcs.
With each arc $j\in\mathcal{A}$, we associate one player, which we
call the $j$th player. The variable controlled by the $j$th player
is the nonnegative integer flow on arc $j$, denoted by $x_{j}\in\mathbf{Z}$.
Each player is trying to minimize a proper cost function $f_{j}:\mathbf{Z}\to\mathbf{R}$,
subject to the network flow constraints. We assume each of the cost
functions is proper, }\textcolor{black}{\emph{i.e.},}\textcolor{black}{{}
for all $i\in[n]$ we have $-\infty\not\not\in f_{i}(\mathbf{Z}),$
and $\textrm{dom}f_{i}=\{x_{i}\in\mathbf{Z}^{n}\mid f_{i}(x_{i})<+\infty\}\neq\emptyset.$
There is an upper bound $u_{j}$, which limits how much flow the $j$th
player can carry through arc $j$. Without any loss of generality,
we take the lower bound on every arc to be $0$ \citep[ Page 39]{Ahuja1988}.
The supply or demand of flow at each node $i\in\mathcal{M}$ is denoted
by $b_{i}$. If $b_{i}>0$, then $i$ is a supply node; if $b_{i}<0$,
then $i$ is a demand node with a demand of $-b_{i}$, and if $b_{i}=0$,
then $i$ is a trans-shipment node. We allow parallel arcs to exist
between two nodes.}
%  The vector formed by all the decision variables
%is denoted by $x=(x_{1},x_{2},\ldots,x_{n})\in\mathbf{Z}^{n}$. By
%$x_{-j}\in\mathbf{Z}^{n-1}$ we denote the vector formed by all the
%players decision variables except $j$th player's decision variable.
%To put emphasis on the $j$th player's variable we sometimes write
%$x$ as $(x_{j},x_{-j})$. Each player has a cost function $f_{j}(x_{j}):\mathbf{Z}\to\mathbf{R}$,
%which depends on its variable $x_{j}$. The goal of the $j$th player
%for $j\in[n]$, given other players' strategies $x_{-j}\in\mathbf{Z}^{n-1}$,
%is to solve its minimization problem subject to the network constraints.

\paragraph*{\textcolor{black}{}}

\textcolor{black}{In any minimum cost flow problem there are three
types of constraints. }

\textcolor{black}{\emph{(i) Mass balance constraint.}}\textcolor{black}{{}
The mass balance constraint states that for any node, the outflow
minus inflow must equal the supply/demand of the node. We describe
the constraint using the \emph{node-arc incidence matrix}. Let us fix a particular
ordering of the arcs, and let $x\in\mathbf{Z}^{n}$ be the resultant vector
of flows.
First, we define the\emph{ augmented node-arc incidence matrix} $\tilde{A}$, where each row corresponds to 
a node, and each column corresponds to an arc. The symbol $\tilde{a}_{ij}$ denotes the $(i,j)$th entry of $\tilde{A}$ that corresponds to the $i$th node and the $j$th arc; $\tilde{a}_{ij}$  is $1$ if $i$
is the start node of the $j$th arc, $-1$ if $i$ is the end node
of the $j$th arc, and $0$ otherwise. Note that parallel arcs will
correspond to different columns with same entries in the augmented node-arc
incidence matrix. So every column of $\tilde{A}$ has exactly two
nonzero entries, one equal to $1$ and one equal to $-1$ indicating
the start node and the end node of the associated arc. Denote, $\tilde{b}=(b_{1},\ldots,b_{m},b_{m+1})$.
Then in matrix notation, we write the mass balance constraint as:
$\tilde{A}x=\tilde{b}.$ The sum of the rows of $\tilde{A}$ is equal to zero
vector, so one of the constraints
associated with the rows of the linear system $\tilde{A}x-\tilde{b}$
is redundant, and by removing the last row of the linear system, we
can arrive at a system, $Ax=b,$ where $A=\tilde{A}_{[1:m,1:n]}$ is the \emph{node-arc incidence matrix},
and $b=\tilde{b}_{[1:m]}$. The vector $b$ is also called the \emph{resource vector}. Now $A$ is a full row rank matrix under
the assumption of $G$ being connected and $\sum_{i\in\mathcal{N}}b_{i}=0$
\citep[Corollary 7.1]{Bertsimas1997}.}\textcolor{black}{\emph{ }}

\textcolor{black}{\emph{(ii) Flow bound constraint.}}\textcolor{black}{{}
The flow on any arc must satisfy the lower bound and capacity constraints,
}\textcolor{black}{\emph{i.e.},}\textcolor{black}{{} $0\preceq x\preceq u$.
The flow bound constraint can often be relaxed or omitted in practice
}\citep[pages 550-551]{Boyd2004}. In such cases, the flow direction is
flexible, and overflow is allowed subject to a suitable penalty.

\textcolor{black}{\emph{(iii) Integrality constraint.}}\textcolor{black}{{}
The flow on any arc is integer-valued, }\textcolor{black}{\emph{i.e.}, $x\in\mathbf{Z}^{n}$.}\textcolor{black}{{}
This }does not incur a significant loss of generality (see Remark
\ref{rem:Integral-arbitrary} below).

\textcolor{black}{So the constraint set, which we denote by $P$ can
be written as, 
\begin{align}
P= & \left\{ x\in\mathbf{Z}^{n}\mid Ax=b,0\preceq x\preceq u\right\} ,\label{eq:constraint-set-P}
\end{align}
}and the subset of $P$ containing only the equality constraints is
denoted by $Q$, \emph{i.e.},
\begin{equation}
Q=\left\{ x\in\mathbf{Z}^{n}\mid Ax=b\right\} .\label{eq:equality-constraint-set-Q}
\end{equation}

\paragraph{}

\textcolor{black}{Consider a set of players denoted by $[n]$. The
decision variable controlled by the $i$th player is $x_{i}\in\mathbf{Z}$,
}\textcolor{black}{\emph{i.e.}},\textcolor{black}{{} each player has
to take an integer-valued action. The vector formed by all the decision
variables is denoted by $x=(x_{1},x_{2},\ldots,x_{n})\in\mathbf{Z}^{n}$.
By $x_{-i}\in\mathbf{Z}^{n-1}$, we denote the vector formed by all
the players decision variables except $i$th player's decision variable.
To put emphasis on the $i$th player's variable we sometimes write
$x$ as $(x_{i},x_{-i})$. Each player has a cost function $f_{i}(x_{i}):\mathbf{Z}\to\mathbf{R}$,
which depends on its variable $x_{i}$. The goal of the $i$th player
for $i\in[n]$, given other players' strategies $x_{-i}\in\mathbf{Z}^{n-1}$,
is to solve the minimization problem}

\textcolor{black}{
\begin{align}
 & \begin{aligned} & \text{minimize}_{x_{i}} &  & f_{i}\left(x_{i}\right)\\
 & \text{subject to} &  & A_{i}x_{i}+\sum_{j=1,j\neq i}^{n}A_{j}x_{j}=b\\
 &  &  & 0\preceq(x_{i},x_{-i})\preceq u\\
 &  &  & x\in\mathbf{Z}^{n}.
\end{aligned}
\label{eq:game_in_consideration}
\end{align}
Our objective is to calculate efficient Pareto optimal points for
the problem. We define vector optimal points first, then Pareto optimal
point, and finally efficient Pareto optimal point.}
\begin{definition}
(\textbf{Vector optimal point}) \textcolor{black}{In problem }\eqref{eq:game_in_consideration}\textcolor{black}{,
a point $x^{\textrm{vo}}\in P$ is }\textcolor{black}{\emph{vector
optimal}}\textcolor{black}{{} if it satisfies the following: }
\[
\left(\forall\tilde{x}\in P\right)\left(\forall i\in[n]\right)\quad f_{i}(\tilde{x}_{i})\geq f_{i}(x_{i}^{\textrm{vo}}).
\]
\end{definition}
\begin{definition}
\textbf{\textcolor{black}{(Pareto optimal point)}}\textcolor{black}{{}
In problem }\eqref{eq:game_in_consideration}\textcolor{black}{, a
point $x^{\textrm{po}}\in P$ is }\textcolor{black}{\emph{Pareto optimal}}\textcolor{black}{{}
if it satisfies the following: there }\textcolor{black}{\emph{does
not }}\textcolor{black}{exist another point $\tilde{x}\in P$ such
that 
\begin{equation}
\left(\forall i\in[n]\right)\quad f_{i}(\tilde{x}_{i})\leq f_{i}(x_{i}^{\textrm{po}}),\label{eq:pareto_def}
\end{equation}
}with at least one index $j\in[n]$ satisfying $f_{j}(\tilde{x}_{j})<f_{j}(x_{j}^{\textrm{po}})$. 
\end{definition}
\begin{definition}
\textbf{\textcolor{black}{(Efficient Pareto optimal point)}} \textcolor{black}{\label{def:Pareto-optimal-point-1}
In problem }\eqref{eq:game_in_consideration}\textcolor{black}{, a
point} $x^{*}$ is an \emph{efficient Pareto optimal solution}, if
it is Pareto optimal and it achieves partial vector optimality over
a maximal subset of $[n]$; \emph{i.e.}, $x^{*}$ satisfies
\eqref{eq:pareto_def} and the set $\mathcal{S}\subseteq[n]$ that
satisfies
\[
\left(\forall\tilde{x}\in P\right)\left(\forall i\in\mathcal{S}\right)\quad f_{i}(\tilde{x}_{i})\geq f_{i}(x_{i}^{*}),
\]
 is maximal.
\end{definition}
\begin{remark}
\label{rem:Integral-arbitrary}In our model, we have taken the flow
through any arc of the network to be an integer. However this assumption
does not incur a significant loss of generality because we can use
our integer model to obtain a real valued Pareto optimal solution
to an arbitrary degree of accuracy by using the following scaling
technique \citep[page 545]{Ahuja1988}. Suppose we want a real valued
Pareto optimal solution $x^{*}$\textcolor{black}{. Such a real valued
Pareto optimal solution corresponds to a modified version of problem
}\eqref{eq:game_in_consideration}\textcolor{black}{{} with the last
constraint being changed to $x\in\mathbf{R}^{n}.$ In practice, we
always have an estimate of how many points after the decimal point
we need to consider. So in the modified problem we substitute each
$x_{i}$ for $i\in[n]$ with $y_{i}/\alpha$, where $y_{i}\in\mathbf{Z}$,
and $\alpha$ is chosen depending on the desired degree of accuracy
(}\textcolor{black}{\emph{e.g.},}\textcolor{black}{{} $\alpha$=$1,000$
or $10,000$ or larger depending on how many points after the decimal
point we are interested in). Then we proceed with our methodology
described in the subsequent sections to compute Pareto optimal solutions
over integers. Let $y^{*}$ be one such integer-valued Pareto optimal
solution. Then $x^{*}=\left(x_{i}^{*}\right)_{i=1}^{n}=\left(\frac{1}{\alpha}y_{i}^{*}\right)_{i=1}^{n}$
corresponds to a real-valued Pareto optimal solution to the degree
of accuracy of $1/\alpha$.}
\end{remark}

\begin{remark}
We can formulate our problem as an $n$ person finite static game
in normal form \citep[pages 88-91]{Basar1995}. In
problem \eqref{eq:game_in_consideration}, a player $i\in[n]$
has a finite but possibly astronomical number of alternatives to choose from the feasible set. Let $\mathfrak{m}_{i}$
be the number of feasible alternatives available to player $i$. Furthermore, define the index set $\mathfrak{M}_{i}=[\mathfrak{m}_{i}]=\{1,\ldots,\mathfrak{m}_{i}\}$
with a typical element of the set designated as $\mathfrak{n}_{i}$,
which corresponds to some flow $x_{i}$. If player $1$ chooses a
strategy $\mathfrak{n}_{1}\in\mathfrak{M}_{1}$, player $2$ chooses a
strategy $\mathfrak{n}_{2}\in\mathfrak{M}_{2}$, and so on for all the other players,
 then the cost incurred to player $i$ is a single number
$\text{\ensuremath{\mathfrak{a}}}_{\mathfrak{n}_{1}\cdots\mathfrak{n}_{n}}^{i}$
that can be determined from problem \eqref{eq:game_in_consideration}.
The ordered tuple of all these numbers (over $i\in[n]$), \emph{i.e.},
$\left(\text{\ensuremath{\mathfrak{a}}}_{\mathfrak{n}_{1}\cdots\mathfrak{n}_{n}}^{1},\text{\ensuremath{\mathfrak{a}}}_{\mathfrak{n}_{1}\cdots\mathfrak{n}_{n}}^{2},\ldots,\text{\ensuremath{\mathfrak{a}}}_{\mathfrak{n}_{1}\cdots\mathfrak{n}_{n}}^{n}\right)$,
constitutes the corresponding unique outcome of the game. For a strategy
$\left(\mathfrak{n}_{1}\cdots\mathfrak{n}_{n}\right)$ that violates
any of the constraints in problem \eqref{eq:game_in_consideration},
the cost is taken as $+\infty$. Players make their decisions independently, and each player unilaterally seeks the minimum possible loss, of course
by also taking into account the possible rational and feasible choices
of the other players. The noncooperative Nash equilibrium solution
concept within the context of this $n$-person game can be described
as follows.
\end{remark}
\begin{definition}
\textbf{(Noncooperative Nash equilibrium)} \citep[page 88]{Basar1995}
\label{def:Noncooperative-Nash-equilibrium} An $n$-tuple of strategies
$\left(\mathfrak{n}_{1}^{{\rm Nash}},\ldots,\mathfrak{n}_{n}^{{\rm Nash}}\right)$
with $\mathfrak{n}_{i}^{{\rm Nash}}\in\mathfrak{M}_{i}$ for all $i\in[n]$,
is said to constitute a noncooperative Nash equilibrium solution for
the aforementioned $n$-person nonzero-sum static finite game in normal form if the following $n$ inequalities are satisfied for all $\mathfrak{n}_{i}\in\mathfrak{M}_{i}$
and all $i\in[n]$: 
\begin{eqnarray}
\mathfrak{a}^{i,{\rm Nash}}=\text{\ensuremath{\mathfrak{a}}}_{\mathfrak{n}_{1}^{{\rm Nash}}\mathfrak{n}_{2}^{{\rm Nash}}\cdots\mathfrak{n}_{n}^{{\rm Nash}}}^{i}\leq\text{\ensuremath{\mathfrak{a}}}_{\mathfrak{n}_{1}^{{\rm Nash}}\mathfrak{n}_{2}^{{\rm Nash}}\cdots\mathfrak{n}_{i}\cdots\mathfrak{n}_{n}^{{\rm Nash}}}^{i}.\label{eq:Nash equilibrium}
\end{eqnarray}
The flow corresponding to $\left(\mathfrak{n}_{1}^{{\rm Nash}},\ldots,\mathfrak{n}_{n}^{{\rm Nash}}\right)$
is denoted by $x^{{\rm Nash}}=\left(x_{1}^{{\rm Nash}},\ldots,x_{n}^{{\rm Nash}}\right)$
and is called the \emph{noncooperative Nash equilibrium flow}. Here, the
$n$-tuple $(\mathfrak{a}^{1,{\rm Nash}},\mathfrak{a}^{2,{\rm Nash}},\ldots,\mathfrak{a}^{n,{\rm Nash}})$
is known as a noncooperative (Nash) equilibrium outcome of the $n$-person
game in normal form. Note that the strategy associated with an efficient
Pareto optimal solution $x^{*}$ in Definition \ref{def:Pareto-optimal-point-1}
also satisfies \eqref{eq:Nash equilibrium} in Definition
\ref{def:Noncooperative-Nash-equilibrium}, thus it is a noncooperative
Nash equilibrium flow. 

\end{definition}

\paragraph*{}

We now extend the class of problems that we are going to investigate,
which is strictly larger than the class defined by problem \eqref{eq:game_in_consideration}
and contains it. We will develop our algorithms for this larger class.
Everything defined in the previous subsection still holds, and we extend the constraint set $P$ and the equality constraint set
$Q$ as follows. The structure of $P$ is still that of a standard
form integer polyhedron, \emph{i.e.}, $P=\{x\in\mathbf{Z}^{n}\mid Ax=b,0\preceq x\preceq u\}$,
where $A$ is a full row rank matrix, but it may not necessarily be a node-arc
incidence matrix only. Denote the convex hull of the points in $P$ by
$\mathop{\bf conv}P$. Consider the relaxed polyhedron $\mathop{\bf relaxed}P=\{x\in\mathbf{R}^{n}\mid Ax=b,0\preceq x\preceq u\}$,
where we have relaxed the condition of $x$ being an integer-valued vector.
We now impose the following assumption.

\begin{assumption} \label{assum:shared_vertex} For any integer-valued $b$,
$\mathop{\bf relaxed}P$ has at least one integer-valued basic solution.\end{assumption} 

As vertices of a polyhedron are also basic solutions \citep[page 50]{Bertsimas1997},
if $\mathop{\bf conv}P$ and $\mathop{\bf relaxed}P$ share at least
one vertex, Assumption \ref{assum:shared_vertex} will be satisfied.
We can see immediately that if $A$ is a node-arc incidence matrix,
then $P$ will belong to this class as $\mathop{\bf conv}P=\mathop{\bf relaxed}P$
for network flow problems \citep[Chapter 19]{Schrijver1998}. In other
practical cases of interest, the matrix can satisfy Assumption \ref{assum:shared_vertex},
\emph{e.g.}, matrices with upper or lower triangular square submatrices
with diagonal entries 1, sparse matrices with $m$ variables appearing
only once with coefficients one, \emph{etc}. Moreover, at the expense
of adding slack variables (thus making a larger dimensional problem),
we can turn the problem under consideration into one satisfying Assumption
\ref{assum:shared_vertex}, though the computational price may be
heavy.

In the rest of the paper, whenever we mention \eqref{eq:constraint-set-P},
\eqref{eq:equality-constraint-set-Q}, and \eqref{eq:game_in_consideration},
they correspond to this larger class of problems containing the network
flow problems, and the full row rank matrix $A$ is associated with
this larger class. So the results developed in the subsequent sections
will hold for a network flow setup.
\begin{remark}
Before proceeding any further, we recall that integer programming problems are $\mathcal{NP}$-hard, and even determining
the existence of one feasible point in $P$ is $\mathcal{NP}$ hard
\citep[page 242]{Bertsimas2005}. So problem \eqref{eq:game_in_consideration}
is at least as hard.
\end{remark}

\section{Problem transformation \label{sec:transforming_the_game}}

In this section, we describe how to transform problem \eqref{eq:game_in_consideration}
into $n-m$ decoupled optimization problems for the last $n-m$ players
and how to reformulate the optimization problems for the rest of the
players using consensus constraints. These transformation and reformulation
are necessary for the development of our algorithms.

\subsection{Decoupling optimization problems for the last $n-m$ players\label{decouple_28th_Dec} }

First, we present the following lemma. Recall that, an integer square
matrix is unimodular if its determinant is $\pm1$.
\begin{lemma} \label{28Dec2018}
Assumption \ref{assum:shared_vertex} holds if and only if we can
extract a unimodular basis matrix from $A$.
\end{lemma}
Without any loss of generality, we rearrange the columns of the matrix
$A$ so that the unimodular basis matrix constitutes the first $m$
columns, \emph{i.e.}, if $A=[A_{1}\mid A_{2}\mid\ldots\mid A_{m}\mid A_{m+1}\mid\ldots\mid A_{n}]$,
then $\det([A_{1}\mid A_{2}\mid\ldots\mid A_{m}])=\pm1$, and we reindex
the variables accordingly. Let us denote $[A_{1}\mid A_{2}\mid\ldots\mid A_{m}]=B$,
so $A=[B\mid A_{m+1}\mid\ldots\mid A_{n}]$. Next we present the following
Lemma.
\begin{lemma}
\label{lem:Equivalence-of-constraint-set}Let, $C=B^{-1}A$ and $d=B^{-1}b$.
Then $C\in\mathbf{Z}^{m\times n}$, $d\in\mathbf{Z}^{m}$, and the sets $Q$ and $P$ (defined in \eqref{eq:constraint-set-P}
and \eqref{eq:equality-constraint-set-Q}, respectively) have the equivalent
representations:
\begin{align}
Q & =\left\{ x\in\mathbf{Z}^{n}\mid Cx=d\right\} ,\label{eq:Q_constraint_set}\\
P & =\left\{ x\in\mathbf{Z}^{n}\mid Cx=d,0\preceq x\preceq u\right\} .\label{eq:constraint_set_v2}
\end{align}
\end{lemma}
Before we present the next result, we recall the following definitions
and facts. A full row rank matrix $A\in\mathbf{Z}^{m\times n}$ 
is in \emph{Hermite normal form}, if it is of the structure $[B|0]$,
where $B\in\mathbf{Z}^{m\times m}$ is invertible and lower triangular,
and $0$ represents $\{0\}^{m \times n-m}$
\citep[Section 4.1]{Schrijver1998}.%
\begin{comment}
The following operations on a matrix are called \emph{elementary integer
column operations}: 
\begin{itemize}
\item adding an integer multiple of one column to another column,
\item exchanging two columns, and
\item multiplying a column by -1. 
\end{itemize}
\end{comment}
{} Any full row rank integer matrix can be brought into the Hermite
normal form using elementary integer column operations in polynomial
time in the size of the matrix \citep[Page 243]{Bertsimas2005}.
\begin{lemma}
The matrix $C$, as defined in Lemma \ref{lem:Equivalence-of-constraint-set},
can be brought into the Hermite normal form $[I|0]$ by elementary
integer column operations, more specifically by adding integer multiple
of one column to another column. \label{lem:The-matrix-C-conversion}
\end{lemma}
\begin{comment}
The steps describing the process described in Lemma \ref{lem:The-matrix-C-conversion}
can be summarized as follows.

\begin{algorithm} 
 % \caption*{Converting $C$ to $[I|0]$}
\label{alg:mat_to_Hermite_normal_form}
\begin{algorithmic}[1] 
 \Procedure {Converting $C$ to $[I|0]$}
\For{$i:= 1,2,\ldots, m$} 
\For{$j:= m+1,m+2,\ldots,n$} 
\State $C_j := C_j-C_{i,j} e_i$ 
\EndFor 
\EndFor 
 \EndProcedure
\end{algorithmic} 
\end{algorithm}
\end{comment}

\begin{lemma}
There exists a unimodular matrix $U$ such that $CU=\left[I\mid0\right]$.
\label{thm:unimod_mat_ext}
\end{lemma}
The following theorem is key in transforming the problem into an equivalent
form with $m$ decoupled optimization problems for players $m+1,m+2,\ldots,n$.
\begin{theorem}
\label{thm:transform-x-to-z}The constraint set $Q$ defined in \eqref{eq:Q_constraint_set} is nonempty. Furthermore, any vector $x\in Q$ can
be maximally decomposed in terms of a new variable $z\in\mathbf{Z}^{n-m}$ as follows:
\begin{equation}
x\in Q\Leftrightarrow\exists z\in\mathbf{Z}^{n-m}\textrm{ such that }x=\begin{bmatrix}d_{1}-h_{1}^{T}z\\
d_{2}-h_{2}^{T}z\\
\vdots\\
d_{m}-h_{m}^{T}z\\
z_{1}\\
\vdots\\
z_{n-m}
\end{bmatrix},\label{eq:x-in-z}
\end{equation}
 where $d_{i}$ is the $i$th component of $d=B^{-1}b$, and $h_{i}^{T}\in\mathbf{Z}^{1 \times n-m}$
is the $i$th row of $B^{-1}A_{[1:m,m+1:n]}$. 
\end{theorem}
We can transform our problem using  the new variable $z$. The advantage
of this transformation is that for player $m+1,m+2,\ldots,n$, we
have decoupled optimization problems. By solving these decoupled problems, we can reduce
the constraint set significantly (especially when the number of minimizers
for the constraint set are small). 

From Theorem \ref{thm:transform-x-to-z}, $x_{i}=z_{i-m}$ for $i\in[m+1:n]$.
For problem \eqref{eq:game_in_consideration}, we can write the optimization
problem for any player $m+i$ for $i\in[n-m]$ as follows.
\begin{align}
 & \begin{aligned} & \text{minimize}_{z_{i}} &  & f_{i}(z_{i})\\
 & \text{subject to} &  & 0\leq z_{i}\leq u_{i}\\
 &  &  & z_{i}\in\mathbf{Z}.
\end{aligned}
\label{eq:decoupled-opt}
\end{align}
Each of these optimization problems is a decoupled univariate optimization
problem, which can be easily solved graphically. We can optimize over
real numbers, find the minimizers of the resultant relaxed optimization
problem, determine whether the floor or ceiling of such a
minimizer results in the minimum cost, and pick that as a minimizer
of the original problem. Solving $n-m$ decoupled optimization problems
immediately reduces the constraint set into a much smaller set. Let
us denote the set of different optimal solutions for player $m+i$
for $i\in[n-m]$ as $D_{i}=\left\{ z_{i,1},z_{i,2},\ldots,z_{i,p_{i}}\right\} $
sorted from smaller to larger, where $p_{i}$ is the total number
of minimizers. Define, $D=\bigtimes_{i=1}^{n-m}D_{i}$. Note that
$D\neq\emptyset$ .

\subsection{Consensus reformulation for the first $m$ players\label{sec:Transforming-the-m-players}}

In this section, we transform the optimization problems for
the first $m$ players using \eqref{eq:x-in-z} and consensus constraints.
Consider the optimization problems for the first $m$ players in variable
$z$, which have coupled costs due to \eqref{eq:x-in-z}. We deal
with the issue by introducing \emph{consensus constraints} \citep[Section 5.2]{Parikh2013}.
We provide each player $i\in[m]$ with its own local copy of $z$,
denoted by $z^{(i)}\in\mathbf{Z}^{n-m}$, which acts as its decision
variable. This local copy has to satisfy the following conditions.
\emph{First}, using \eqref{eq:x-in-z} for any $i\in[m]$, $x_{i}=d_{i}-h_{i}^{T}z^{(i)}$.
The copy $z^{(i)}$ has to be in consensus with the rest of the first
$m$ players, \emph{i.e.}, $z^{(i)}=z^{(j)}$ for all $j\in[m]\setminus\{i\}$.\emph{
Second}, the copy $z^{(i)}$ has to satisfy the flow bound constraints,
\emph{i.e.}, $0\leq d_{i}-h_{i}^{T}z^{(i)}\leq u_{i}$ for all $i\in[m]$.
\emph{Third}, for the last $n-m$ players $z_{i}\in D_{i}$, as obtained
from the solutions of the decoupled optimization problems \eqref{eq:decoupled-opt},
so $z^{(i)}$ has to be in $D$, \emph{i.e.}, for all $i\in[m]$ we
have\emph{ 
\[
z^{(i)}=z\in D\Leftrightarrow\left(\forall j\in[n-m]\right)z_{j}^{(i)}\in D_{j}.
\]
 }Combining the aforementioned conditions, for all $i\in[m]$, the $i$th player's optimization problem
in variable $z^{(i)}$ can be written as: 
\begin{align}
 & \begin{aligned} & \text{minimize}_{z^{(i)}} &  & \bar{f}_{i}\left(z^{(i)}\right)=f_{i}(d_{i}-h_{i}^{T}z^{(i)})\\
 & \text{subject to} &  & z^{(i)}=z^{(j)},\quad j\in[m]\setminus\{i\}\\
 &  &  & 0\leq d_{i}-h_{i}^{T}z^{(i)}\leq u_{i}\\
 &  &  & z_{j}^{(i)}\in D_{j},\quad j\in[n-m].
\end{aligned}
\label{eq:intermediate-in-m}
\end{align}
An integer linear inequality constraint $\alpha\leq v\leq\beta,$
where $\alpha,\beta,v\in\mathbf{Z}$, is equivalent to $v\in\{\alpha,\alpha+1,\ldots,\beta\}\Leftrightarrow(v-\alpha)(v-\alpha-1)\cdots(v-\beta)=0$.
Using this fact, we write the last two constraints in \eqref{eq:intermediate-in-m}
in polynomial forms as follows.
\begin{align}
 & q_{i}(z^{(i)})=(d_{i}-h_{i}^{T}z^{(i)})(d_{i}-h_{i}^{T}z^{(i)}-1)\cdots(d_{i}-h_{i}^{T}z^{(i)}-u_{i})=0,\label{eq:Hilbert_system}\\
 & r_{j}(z^{(i)})=(z_{j}^{(i)}-z_{j,1})(z_{j}^{(i)}-z_{j,2})\ldots(z_{j}^{(i)}-z_{j,p_{i}})=0,\quad\quad j\in[n-m].\label{eq:Hilbert_system_2}
\end{align}
Hence for all $i\in[m]$ any feasible $z^{(i)}$ for problem \eqref{eq:intermediate-in-m}
comes from the following set: 
\begin{eqnarray}
\mathcal{F} & = & \bigcap_{k=1}^{m}\{z\in\mathbf{Z}^{n-m}\mid q_{k}(z)=0, \textrm{and} \left(\forall j\in[n-m]\right)\quad r_{j}(z)=0\}\nonumber \\
 & = & \{z\in\mathbf{Z}^{n-m}\mid\left(\forall k\in[m]\right)\;q_{k}(z)=0, \textrm{and} \left(\forall j\in[n-m]\right)\;r_{j}(z)=0\}.\label{eq:mathcal-F}
\end{eqnarray}
In \eqref{eq:mathcal-F}, the intersection in the first line ensures
that the consensus constraints are satisfied, and the second line
just expands the first. So the optimization problem \eqref{eq:intermediate-in-m}
is equivalent to 
\begin{align}
 & \begin{aligned} & \textrm{minimize}_{z^{(i)}} &  & \bar{f}_{i}\left(z^{(i)}\right)\\
 & \text{subject to} &  & z^{(i)}\in\mathcal{F},
\end{aligned}
\label{eq:final-optimization-model-3}
\end{align}
for $i\in[m]$. Thus each of these players is optimizing
over a \emph{common constraint set} $\mathcal{F}$. So finding the
points in $\mathcal{F}$ is of interest, which we discuss next. 

\section{Algorithms \label{sec:Algorithms-to-compute}}

In this section, first, we review some necessary background on algebraic
geometry, and then we present a theorem to check if $\mathcal{F}$
is nonempty and provide an algorithm to compute the points in a nonempty
$\mathcal{F}$. Finally, we present our algorithm to compute efficient
Pareto optimal points. In devising our algorithm, we use algebraic
geometry rather than integer programming techniques for the following
reasons. \emph{First,} in this way, we are able to provide an algebraic
geometric characterization for the set of efficient Pareto optimal
solutions for our problem. \emph{Second,} we can show that this set
is nonempty if and only if\emph{ the reduced Groebner basis }(disucussed
in Section \ref{sec:Background} below) of a certain set associated
with the problem is not equal to $\{1\}$ (Theorem \ref{lem:poly_sys_ideal_lemma}).
\emph{Third, }the mentioned result has an algorithmic significance:
the reduced Groebner basis can be used to construct algorithms to
calculate efficient Pareto optimal points.

\subsection{\label{sec:Background}Background on algebraic geometry}

A \emph{monomial} in variables $x=(x_{1},x_{2},\ldots,x_{n})$ is
a product of the structure $x^{\alpha}=x_{1}^{\alpha_{1}}\cdots x_{n}^{\alpha_{n}},$
 where $\alpha=(\alpha_{1},\ldots,\alpha_{n})\in\mathbf{N}{}^{n}$.
A \emph{polynomial} is an expression that is the sum of a finite number
of terms with each term being a monomial times a real or complex
coefficient. The set of all real polynomials in $x=(x_{1},\ldots,x_{n})$
with real and complex coefficients are denoted by $\mathbf{R}[x]$
and $\mathbf{C}[x]$, respectively, with the variable ordering $x_{1}>x_{2}>\cdots>x_{n}$.
The ideal generated by $f_{1},f_{2},\ldots,f_{m}\in\mathbf{\mathbf{C}}[x]$
is the set 
\[
\mathbf{ideal}\left\{ f_{1},\ldots,f_{m}\right\} =\left\{ \sum_{i=1}^{m}h_{i}f_{i}\mid\left(\forall i\in[m]\right)\;h_{i}\in\mathbf{C}[x]\right\} .
\]
Consider $f_{1},f_{2},\ldots,f_{s}$ which are polynomials in $\mathbf{C}[x]$.
The \emph{affine variety} $V$ of $f_{1},f_{2},\ldots,f_{m}$ is given
by 
\begin{eqnarray}
V(f_{1},\ldots,f_{m}) & = & \left\{ x\in\mathbf{C}^{n}\mid\left(\forall i\in[m]\right)\quad f_{i}(x)=0\right\} .\label{eq:aff_variety_def}
\end{eqnarray}

A \emph{monomial order} on $\mathbf{C}[x_{1},\ldots,x_{n}]$ is a
relation, denoted by $\succ$, on the set of monomials $x^{\alpha},\alpha\in\mathbf{N}^{n}$
satisfying the following. \emph{First}, it is a total order; \emph{second}, if $x^{\alpha}\succ x^{\beta}$
and $x^{\gamma}$ is any monomial, then $x^{\alpha+\gamma}\succ x^{\beta+\gamma}$; \emph{third}, every nonempty subset of $\mathbf{N}^{n}$ has a smallest element
under $\succ$. We will use \emph{lexicographic order}, where we say
$x^{\alpha}\succ_{\text{lex}}x^{\beta}$ if and only if the left most
nonzero entry of $\alpha-\beta$ is positive. Suppose that we are given
a monomial order $\succ$ and a polynomial $f(x)=\sum_{\alpha\in S}f_{\alpha}x^{\alpha}$.
The \emph{leading term} of the polynomial with respect to $\succ$,
denoted by $\textrm{lt}_{\succ}\left(f\right)$, is that monomial
$f_{\alpha}x^{\alpha}$ with $f_{\alpha}\neq0$, such that $x^{\alpha}\succ x^{\beta}$
for all other monomials $x^{\beta}$ with $f_{\beta}\neq0$. The monomial
$x^{\alpha}$ is called the \emph{leading monomial} of $f$. Consider
a nonzero ideal $I$. The set of the leading terms for the polynomials
in $I$ is denoted by $\textrm{lt}_{\succ}\left(I\right)$. Thus $\textrm{lt}_{\succ}\left(I\right)=\left\{ cx^{\alpha}\mid\left(\exists f\in I\right)\quad\textrm{lt}_{\succ}(f)=cx^{\alpha}\right\} .$
By $\mathbf{ideal}\left\{ \textrm{lt}_{\succ}\left(I\right)\right\} $
with respect to $\succ$, we denote the ideal generated by the elements
of $\textrm{lt}_{\succ}\left(I\right)$. 

A \emph{Groebner basis} $G_{\succ}$ of an ideal $I$ with respect
to the monomial order $\succ$ is a finite set of polynomials $g_{1},\ldots,g_{t}\in I$
such that $\mathbf{ideal}\left\{ \textrm{lt}_{\succ}\left(I\right)\right\} =\mathbf{ideal}\left\{ \textrm{lt}_{\succ}\left(g_{1}\right),\ldots,\textrm{lt}_{\succ}\left(g_{t}\right)\right\} .$
A \emph{reduced Groebner basis} $G_{\textrm{reduced},\succ}$ for
an ideal $I$ is a Groebner basis for $I$ with respect to monomial
order $\succ$ such that, for any $f\in G_{\textrm{reduced},\succ}$,
the coefficient associated with $\textrm{lt}_{\succ}\left(f\right)$
is $1$, and for all $f\in G_{\textrm{reduced},\succ}$, no monomial
of $f$ lies in $\mathbf{ideal}\left\{ \textrm{lt}_{\succ}\left(G\setminus\{f\}\right)\right\} $.
For a nonzero ideal $I$ and given monomial ordering, the reduced Groebner
basis is unique \citep[Proposition 6,][Page 92]{Cox2007}. Suppose
$I=\mathbf{ideal}\left\{ f_{1},\ldots,f_{m}\right\} \subseteq\mathbf{C}\left[x_{1},\ldots,x_{n}\right]$.
Then for any $l\in[n]$, the $l$th \emph{elimination ideal} for $I$
is defined by $I_{l}=I\cap\mathbf{C}[x_{l+1},\ldots,x_{n}].$ Let $I\subseteq\mathbf{C}\left[x_{1},\ldots,x_{n}\right]$
be an ideal, and let $G$ be a Groebner basis of $I$ with respect
to lexicographic order with $x_{1}\succ x_{2}\succ\ldots\succ x_{n}$.
Then for every integer $l\in\left\{ 0,n-1\right\}$, the set $G_{l}=G\cap\mathbf{C}\left[x_{l+1},\ldots,x_{n}\right]$
is a Groebner basis for the $l$th elimination ideal $I_{l}$.

\subsection{Nonemptyness of $\mathcal{F}$ \label{sec:Nonemptyness-of-}}

We will use the following theorem in proving the results in this section.
\begin{theorem}
(Weak Nullstellensatz \citep[Theorem 1,][page 170]{Cox2007}) Consider
$f_{1},f_{2},\ldots,f_{s}$ as polynomials in $\mathbf{C}[x_{1},x_{2},\ldots,x_{n}]$.
If $V(f_{1},f_{2},\ldots,f_{m})=\emptyset$ (see \eqref{eq:aff_variety_def}
for definition of $V$), then 
\[
\mathbf{ideal}\{f_{1},f_{2},\ldots,f_{m}\}=\mathbf{C}[x_{1},x_{2},\ldots,x_{n}].
\]
\end{theorem}
First, we present the following result.
\begin{theorem}
\label{lem:poly_sys_ideal_lemma} The set $\mathcal{F}$ is nonempty
if and only if
\[
G_{\textrm{reduced},\succ}\ne\{1\},
\]
 where $G_{\textrm{reduced},\succ}$ is the reduced Groebner basis
of $\mathbf{ideal}\left\{ q_{1},\ldots,q_{m},r_{1},\ldots,r_{n-m}\right\} $
with respect to any ordering.
\end{theorem}
\begin{remark}
\label{rem:cplx_eq_int}In the proof, we have shown that
feasibility of the system \eqref{eq:poly_sys_2-2} in $\mathbf{C}^{n-m}$
is equivalent to its feasibility in $\mathbf{Z}^{n-m}$.
So
\begin{equation}
V\left(q_{1},\ldots,q_{m},r_{1},\ldots,r_{n-m}\right)\cap\mathbf{Z}^{n-m}=V\left(q_{1},\ldots,q_{m},r_{1},\ldots,r_{n-m}\right).\label{eq:variety equality}
\end{equation}

If we are interested in just verifying the feasibility of the polynomial
system, then calculating a reduced Groebner basis with respect to
any ordering suffices. However, if we are interested in extracting
the feasible points, then we choose lexicographic ordering as lexicographic
ordering allows us to use algebraic elimination theory. There are
many computer algebra packages to compute reduced Groebner basis such
as Macaulry2, SINGULAR, FGb, Maple, and Mathematica. We now
describe how to extract the points in $\mathcal{F}$.
\end{remark}

Suppose $G_{\textrm{reduced},\succ}\ne\{1\}$. Naturally, the next
question is how to compute points in $\mathcal{F}$? In the next section,
we will show that the points in $\mathcal{F}$ are related to the
efficient Pareto optimal points that we are seeking. We now briefly
discuss systematic methods for extracting $\mathcal{F}$ based on
algebraic elimination theory, a branch of computational algebraic
geometry. For details on elimination theory, we refer the interested
readers to \citep[Chapter 3]{Cox2007}. First, we present the following
lemma.

\begin{lemma} \label{choke_amar_trishna}
Suppose $G_{\textrm{reduced},\succ}\ne\{1\}$. Then $\mathcal{F}=V\left(G_{\textrm{reduced},\succ_{\textrm{lex}}}\right)\neq\emptyset.$
\end{lemma}

Algorithm \ref{alg:extracting_mathcal_F-2} below calculates all the
points in $\mathcal{F}$, when $G_{\textrm{reduced},\succ_{\textrm{lex}}}\neq\{1\}.$
Lemma \ref{lem:Algorithm--correctly} proves its accuracy.

\begin{lemma}
Algorithm \ref{alg:extracting_mathcal_F-2} correctly calculates all
the points in $\mathcal{F}$, when it is nonempty. \label{lem:Algorithm--correctly}
\end{lemma}

\textcolor{black}{\scriptsize{}}
\begin{algorithm}
\textbf{\textcolor{black}{Input:}}\textcolor{black}{{} Polynomial system
$q_{i}(z)=0$ for $i\in[m]$, and $r_{j}(z)=0$ for $j\in[n-m]$, $G_{\textrm{reduced},\succ_{\textrm{lex}}}\neq\{1\}.$}

\textbf{\textcolor{black}{Output:}}\textcolor{black}{{} The set $\mathcal{F}$.}

\textcolor{black}{\hrulefill}
\begin{description}
\item [{Step}] 1.
\begin{itemize}
\item \textcolor{black}{Calculate the set $G_{n-m-1}=G_{\textrm{reduced},\succ_{\textrm{lex}}}\cap\mathbf{C}\left[z_{n-m}\right],$ which
is a Groebner basis of the $(n-m)$th elimination ideal of $\mathbf{ideal}\left\{ q_{1},\ldots,q_{m},r_{1},\ldots,r_{n-m}\right\} $ that  consists of univariate polynomials in $z_{n-m}$ as an implication
of \citep[page 116, Theorem 2]{Cox2007}. }
\item \textcolor{black}{Find the variety of $G_{n-m-1}$, denoted by $V\left(G_{n-m-1}\right),$
which contains the list all possible $z_{n-m}$ coordinates for
the points in $\mathcal{F}$.}
\end{itemize}
\item [{\textcolor{black}{Step}}] \textcolor{black}{2.}
\begin{itemize}
\item \textcolor{black}{Calculate $G_{n-m-2}=G_{\textrm{reduced},\succ_{\textrm{lex}}}\cap\mathbf{C}\left[z_{n-m-1},z_{n-m}\right],$
which is again a Groebner basis of the $(n-m-1)$th elimination ideal
of $\mathbf{ideal}\left\{ q_{1},\ldots,q_{m},r_{1},\ldots,r_{n-m}\right\} $ and
consists of bivariate polynomials in $z_{n-m}$ and $z_{n-m-1}$. }
\item \textcolor{black}{From Step 1, we already have the $z_{n-m}$ coordinates
for the points in $\mathcal{F}$. So by substituting those $|V(G_{n-m-1})|$
values in $G_{n-m-2}$, we arrive at a set of univariate
polynomials in $z_{n-m-1}$, which we denote by $\{\bar{G}_{n-m-2}^{(i)}\}_{i=1}^{|V(G_{n-m-1})|}$. }
\item \textcolor{black}{For all $i=1,2,\ldots,|V(G_{n-m-1})|$, find the
variety of $\bar{G}_{n-m-2}^{(i)}$, denoted by $V(\bar{G}_{n-m-2}^{(i)})$,
which contains the list all possible $z_{n-m-1}$ coordinates
associated with a particular $z_{n-m}\in V(G_{n-m-1})$. We now have
all the possible $\left(z_{n-m-1},z_{n-m}\right)$ coordinates of
$\mathcal{F}$.}
\end{itemize}
\item [{Step}] 3.
\begin{itemize}
\item \textcolor{black}{We repeat this procedure for $G_{n-m-3},G_{n-m-4},\ldots,G_{0}$.
In the end, we have set of all points in $\mathcal{F}$. }
\end{itemize}
\end{description}
\textbf{\textcolor{black}{return}}\textcolor{black}{{} $\mathcal{F}$.}

\textcolor{black}{\caption{Extracting the points in $\mathcal{F}$}
\label{alg:extracting_mathcal_F-2}}
\end{algorithm}
{\scriptsize \par}

\subsection{Finding efficient Pareto optimal points from $\mathcal{F}$\label{sec:Finding-the-Pareto}}

Suppose $G_{\textrm{reduced},\succ_{\textrm{lex}}}\neq\{1\}$, and
using Algorithm \ref{alg:extracting_mathcal_F-2} we have computed
$\mathcal{F}$. We now propose Algorithm \ref{alg:other-m-players-1}
and show that the resultant points are Pareto optimal.

\begin{algorithm}
\textbf{Input:} The optimization problem \eqref{eq:final-optimization-model-3}
for any $i\in[m]$, $\mathcal{F}\neq\emptyset$.

\textbf{Output:} Efficient Pareto optimal solutions for problem \eqref{eq:game_in_consideration}.

\hrulefill
\begin{description}
\item [{Step}] 1.
\begin{description}
\item [{\textbf{for}}] $i=1,\ldots,m$

\[
X_{i}:=\left\{ d_{i}-h_{i}^{T}z^{(i)}\mid z^{(i)}\in\mathcal{F}\right\} =d_{i}-h_{i}^{T}\mathcal{F},
\]

\begin{equation}
\left(\forall x_{i}\in X_{i}\right)\quad(X_{i})^{-1}(x_{i}):=\{z^{(i)}\in\mathcal{F}\mid x_{i}=d_{i}-h_{i}^{T}z^{(i)}.\}\label{eq:xi-inv}
\end{equation}
 
\item [{\textbf{end}}] \textbf{for}
\end{description}
\item [{Step}] 2.

Sort the elements of the $\{X_{i}\}_{i=1}^{m}$s with respect to cardinality
of the elements in a descending order. Denote the index set of the
sorted set by $\{s_{1},s_{2},\ldots,s_{m}\}$ \texttt{}such that
\[
|X|_{s_{1}}\geq|X|_{s_{2}}\geq\cdots\geq|X|_{s_{m}}.
\]

\item [{Step}] 3.
\begin{description}
\item [{\textbf{\textcolor{black}{for}}}] \textcolor{black}{$i\in[m]$}

\textcolor{black}{Solve the univariate optimization problem 
\begin{align}
 & \begin{aligned} & \mathrm{minimize}_{x_{s_{i}}} &  & f_{s_{i}}(x_{s_{i}})\\
 & \mathrm{subject\;to} &  & x_{s_{i}}\in X_{s_{i}},
\end{aligned}
\label{eq:univar-last-m-players}
\end{align}
 and denote the set of solutions by $X_{s_{i}}^{*}$. Set}

\textcolor{black}{
\begin{eqnarray}
\mathcal{F}_{s_{i}}^{*} & := & \bigcup_{x_{s_{i}}\in X_{s_{i}}^{*}}(X_{s_{i}}^{*})^{-1}(x_{s_{i}})\subseteq\mathcal{F},\label{eq:mcFsi}
\end{eqnarray}
}
\begin{description}
\item [{\textbf{if}}] $i\leq m$

\begin{align}
 & X_{s_{i+1}}:=\left\{ d_{s_{i+1}}-h_{s_{i+1}}^{T}z\mid z\in\mbox{\ensuremath{\mathcal{F}}}_{s_{i}}^{*}\right\} .\label{eq:Xsi1}
\end{align}

\item [{\textbf{end}}] \textbf{if}
\end{description}
\item [{\textbf{end}}] \textbf{for}
\end{description}
\end{description}
\textbf{return} $\mathcal{F}_{s_{m}}^{*}.$

\caption{\label{alg:other-m-players-1} Computing the set of solutions to problem
\eqref{eq:final-optimization-model-3}. }
\end{algorithm}

We have the following results for Algorithm \ref{alg:other-m-players-1}.
\begin{lemma}
\label{lem:subseteq} In Algorithm \ref{alg:other-m-players-1}, for
all $i\in[m-1]$, we have $\mathcal{F}_{s_{i+1}}^{*}\subseteq\mathcal{F}_{s_{i}}^{*}\subseteq\mathcal{F}$
($\mathcal{F}_{s_{i}}^{*}$ defined in \eqref{eq:mcFsi}).
\end{lemma}

\begin{lemma}
\label{lem:optimal-x-z}In Algorithm \ref{alg:other-m-players-1},
for any $i\in[m]$, $x_{s_{i}}\in X_{s_{i}}^{*}$ (the set $X_{i}^{*}$
is defined in Step 1, 3 of Algorithm \ref{alg:other-m-players-1})
if and only if $z^{*}\in\mathcal{F}_{s_{i}}^{*}$. Furthermore, $z^{*}\in\mathcal{F}_{s_{i}}^{*}$
solves the following optimization problem 
\begin{align*}
 & \begin{aligned} & \mathrm{minimize}_{z} &  & f_{s_{i}}(d_{s_{i}}-h_{s_{i}}^{T}z)\\
 & \mathrm{subject\;to} &  & z\in\mathcal{F}_{s_{i-1}}^{*},
\end{aligned}
\end{align*}
for all $i\in[2:m]$. 
\end{lemma}
In Algorithm \ref{alg:other-m-players-1}, at no stage
can $\mathcal{F}_{s_{i}}^{*}$ get empty.
\begin{lemma}
Suppose $\mathcal{F}\neq\emptyset$. Then in Algorithm \ref{alg:other-m-players-1},
$\mathcal{F}_{s_{i}}^{*}$ is nonempty for any $i\in[m]$.
\end{lemma}

\begin{theorem}
For any $z^{*}\in\mathcal{F}_{s_{m}}^{*}$, 
\begin{equation}
x^{*}=(d_{1}-h_{1}^{T}z^{*},\ldots,d_{m}-h_{m}^{T}z^{*},z_{1}^{*},\ldots,z_{n-m}^{*})\label{eq:x*z*}
\end{equation}
 is an efficient Pareto optimal point.
\end{theorem}

\section{A penalty based approach when $\text{\ensuremath{\mathcal{F}}}$
is empty\label{sec:A-penalty-based}}

In the case that $\mathcal{F}$ is empty, we design a penalty based
approach to solve a penalized version of problem \eqref{eq:final-optimization-model-3}
that can be of use to network administrators and policy makers. First,
 for $i\in[m]$, using \eqref{eq:mathcal-F},
\eqref{eq:Hilbert_system_2}, and \eqref{eq:Hilbert_system}, problem
\eqref{eq:final-optimization-model-3} can be written in the following
equivalent form: %

\begin{align}
 & \begin{aligned} & \textrm{minimize}_{z^{(i)}\in D} &  & \bar{f}_{i}\left(z^{(i)}\right)\\
 & \text{subject to} &  & \left(\forall k\in[m]\right)\quad q_{k}(z^{(i)})=0.
\end{aligned}
\label{eq:final-optimization-model-3-1}
\end{align}

When $\mathcal{F}=\emptyset$, we relate the problem above to a penalized
version, which is a standard practice in operations research literature
\citep[Page 278]{opac-b1108032}. In this penalized version, we disregard
the equality constraints $q_{k}(z^{(i)})=0$ for $k\in[m]$, rather
we we augment the cost with a term that penalizes the violation of
these equality constraints. So a penalized version of the problem
\eqref{eq:final-optimization-model-3-1} is as follows:

\begin{align}
 & \begin{aligned} & \textrm{minimize}_{z^{(i)}} &  & \bar{f}_{i}\left(z^{(i)}\right)+\sum_{k=1}^{m}\gamma_{k}p\left(q_{k}(z^{(i)})\right),\\
 & \textrm{subject to} &  & z^{(i)}\in D.
\end{aligned}
\label{eq:final-optimization-model-3-1-1}
\end{align}
where $p:\mathbf{R}\to\mathbf{R}_{+}$ is a penalty function, and
$\gamma_{k}$ is a positive penalty parameter. Some common penalty
functions are:
\begin{itemize}
\item exact penalty: $p:x\mapsto x^{2}$,
\item power barrier: $p:x\mapsto|x|$, \emph{etc.}
\end{itemize}
We have already shown $D$ to be a nonempty set. From a network point
of view, the penalized problem \eqref{eq:final-optimization-model-3-1-1}
has the following interpretation. For the first $m$ arcs in the directed
network under consideration, rather than having a strict \textcolor{black}{\emph{flow
bound constraint}} (see the discussion in Section \ref{subsec:Problem-statement-for-networks}), \textcolor{black} we
have a penalty when $x_{i}<0$ or $x_{i}>u_{i}$ for $i\in[m]$. The
flow bound constraint is still maintained for $i\in[m+1:n]$. In this
regard, the original problem defined by \eqref{eq:game_in_consideration}
has the following penalized version:

\textcolor{black}{
\begin{align}
 & \begin{aligned} & \text{minimize}_{x_{i}} &  & f_{i}\left(x_{i}\right)+\left({\textrm{penalty for violation of }\atop 0\leq x_{i}\leq u_{i}\textrm{ for player }i\in[m]}\right)\\
 & \text{subject to} &  & Ax=A(x_{i},x_{-i})=A_{i}x_{i}+\sum_{j=1,j\neq i}^{n}A_{j}x_{j}=b\\
 &  &  & 0\leq x_{i}\leq u_{i},\quad i\in[m+1:n]\\
 &  &  & x\in\mathbf{Z}^{n}.
\end{aligned}
\label{eq:game_in_consideration-1}
\end{align}
}

Problem \eqref{eq:game_in_consideration-1} can be considered as a
network flow problem, where flow direction is flexible or overflow is allowed subject to a suitable penalty, which is a quite
realistic scenario in practice \citep[pages 550-551]{Boyd2004}. With
this penalized problem, we can proceed as follows. In the developments
of Section \ref{sec:Finding-the-Pareto} set: 
\begin{align*}
 & \bar{f}_{i}\left(z^{(i)}\right):=\bar{f}_{i}\left(z^{(i)}\right)+\sum_{k=1}^{m}\gamma_{k}p\left(q_{k}(z^{(i)})\right),\\
 & \mathcal{F}:=D\neq\emptyset,
\end{align*}
 and then apply Algorithm \ref{alg:other-m-players-1}, which will
calculate efficient Pareto optimal points for the penalized problem
\eqref{eq:game_in_consideration-1}. 

The described penalty scheme can be of use to network administrators
and policy makers to enforce a decision making architecture. Such
an architecture would allow the players to make independent decisions
while ensuring that \emph{(i)} total amount of flow is conserved in
the network by maintaining the \textcolor{black}{\emph{mass balance
constraint}}, \emph{(ii)} the \emph{flow bound constraint} is strictly
enforced for the last $n-m$ players and is softened for the first
$m$ players by imposing penalty, and yet \emph{(iii)} an efficient
Pareto optimal point for the penalized problem can be achieved by
the players \textcolor{black}{where none of their objective functions
can be improved without worsening some of the other players' objective
values. }

\textcolor{black}{From both the players' and the policy maker's point
of views, the penalty based approach makes sense and can be considered
fair for the following reason. In the exact version, each of the last
$n-m$ players gets to minimize its optimization problem }\eqref{eq:decoupled-opt}\textcolor{black}{{}
in a decoupled manner}, whereas each of the first $m$ players is
solving a more restrictive optimization problem \eqref{eq:univar-last-m-players}.
So cutting each of the first $m$ players some slack by softening the flow
bound constraint, where it can carry some extra flow or flow in
opposite direction by paying a penalty, can be considered fair.
\begin{figure}
\centering{}\includegraphics[scale=0.6]{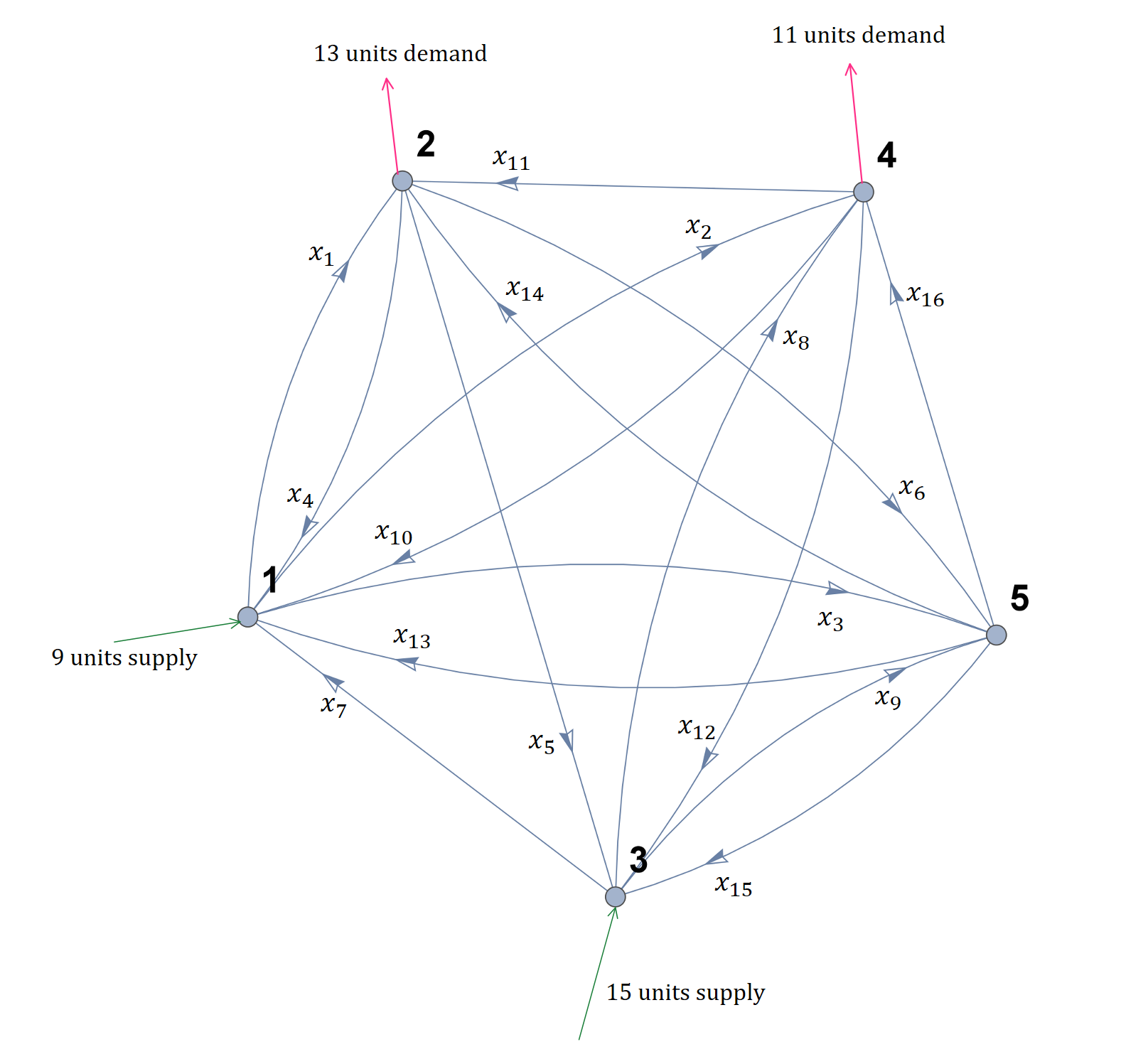}\caption{Network under consideration\label{fig:network-in-consideration}}
\end{figure}
\section{Numerical example\label{sec:Numerical-Examples}}
In this section, we present an illustrative numerical example of our methodology 
in a transportation setup. We have used \texttt{Wolfram Mathematica 10} for numerical computation.
\paragraph*{Problem setup.}$\;$Consider the following directed network with 5 nodes and 16 arcs as
shown in Figure \ref{fig:network-in-consideration}. Nodes 2 and 4
represent two retail centers with demands for 13 and 11 units of a
certain product. The warehouses are denoted by nodes 1 and 3, which
supply 9 and 15 units, respectively. Node 5 is a trans-shipment node.
Different modes of shipment from one node to other is represented
by the arcs in the figure, and these shipments are carried
out by different organizations (carriers). The cost of
a certain shipment depends on the number of products shipped; it is
nonlinear and not necessarily convex. With each arc, we associate one
carrier (player). 

The cost functions for the players are listed in Table \ref{tab:Cost-functions-network}.  The cost functions associated with players 5, 7, 10, and 15 are convex, and the cost functions associated with the rest of the players are nonconvex. Each of the players is trying to minimize its cost. The number of products carried by each player has to be non-negative, and the capacity of each player is represented by 
\[
u=(10,7,11,13,16,12,4,5,6,14,13,15,5,6,6,10),
\]
where $u_i$ represents the maximum capacity for the number of products carried by player $i$.

\begin{table}
\begin{tabular}{|c|c|}
\hline 
Player & Cost function\tabularnewline
\hline 
\hline 
1 & $-\frac{x_{1}^{4}}{30}-\frac{13x_{1}^{3}}{15}+\frac{259x_{1}^{2}}{30}-\frac{263x_{1}}{15}+1$\tabularnewline
\hline 
2 & $\frac{77x_{2}^{5}}{120}-\frac{247x_{2}^{4}}{24}+\frac{471x_{2}^{3}}{8}-\frac{3365x_{2}^{2}}{24}+\frac{6779x_{2}}{60}+1$\tabularnewline
\hline 
3 & $\frac{47x_{3}^{4}}{24}-\frac{133x_{3}^{3}}{4}+\frac{4897x_{3}^{2}}{24}-\frac{2123x_{3}}{4}+485$\tabularnewline
\hline 
4 & $\frac{323x_{4}^{5}}{3360}-\frac{2179x_{4}^{4}}{1120}+\frac{47393x_{4}^{3}}{3360}-\frac{48709x_{4}^{2}}{1120}+\frac{7885x_{4}}{168}+5$\tabularnewline
\hline 
5 & $(x_{5}-1)^{2}$\tabularnewline
\hline 
6 & $-\frac{x_{6}^{4}}{8}+\frac{25x_{6}^{3}}{12}-\frac{71x_{6}^{2}}{8}+\frac{95x_{6}}{12}+10$\tabularnewline
\hline 
7 & $|x_{7}-5|$\tabularnewline
\hline 
8 & $\frac{11x_{8}^{7}}{1260}-\frac{7x_{8}^{6}}{36}+\frac{119x_{8}^{5}}{72}-\frac{479x_{8}^{4}}{72}+\frac{4609x_{8}^{3}}{360}-\frac{803x_{8}^{2}}{72}+\frac{155x_{8}}{28}+1$\tabularnewline
\hline 
9 & $-\frac{15}{16}x_{9}^{3}+\frac{365x_{9}^{2}}{16}-\frac{2865x_{9}}{16}+\frac{7315}{16}$\tabularnewline
\hline 
10 & $(x_{10}-10)^{2}$\tabularnewline
\hline 
11 & $\frac{5x_{11}^{4}}{6}-\frac{35x_{11}^{3}}{3}+\frac{355x_{11}^{2}}{6}-\frac{370x_{11}}{3}+90$\tabularnewline
\hline 
12 & $\frac{5x_{12}^{4}}{6}-\frac{25x_{12}^{3}}{3}+\frac{175x_{12}^{2}}{6}-\frac{110x_{12}}{3}+15$\tabularnewline
\hline 
13 & $\frac{5x_{13}^{4}}{6}-15x_{13}^{3}+\frac{595x_{13}^{2}}{6}-280x_{13}+285$\tabularnewline
\hline 
14 & $\frac{5x_{14}^{4}}{6}-\frac{85x_{14}^{3}}{3}+\frac{2155x_{14}^{2}}{6}-\frac{6020x_{14}}{3}+4165$\tabularnewline
\hline 
15 & $|x_{15}-7|$\tabularnewline
\hline 
16 & $\begin{cases}
x_{16}+1, & \text{if }0\leq x_{16}\le3\\
0, & \textrm{if }4\leq x_{16}\leq6\\
(x_{16}+1)^{3}, & \textrm{if }7\leq x_{16}\leq9\\
-\frac{x_{16}^{3}}{6}+\frac{13x_{16}^{2}}{2}-\frac{244x_{16}}{3}+330, & \textrm{else}
\end{cases}$\tabularnewline
\hline 
\end{tabular}\caption{Cost functions for the players in the network considered \label{tab:Cost-functions-network}}
\end{table}

We seek efficient Pareto optimal points in this setup. 

\paragraph*{Computing node-arc incidence matrix and resource vector.} $\;$First, we compute the node-arc incidence matrix associated with the network under consideration by following the procedure mentioned in the description of the mass balance constraint in Section \ref{subsec:Problem-statement-for-networks}. The resultant augmented node-arc incidence matrix of the network is: 
\[ 
\tilde{A}=   
\begin{pmatrix}
-1 & -1 & -1 & 1 & 0 & 0 & 1 & 0 & 0 & 1 & 0 & 0 & 1 & 0 & 0 & 0 \\ 
1 & 0 & 0 & -1 & -1 & -1 & 0 & 0 & 0 & 0 & 1 & 0 & 0 & 1 & 0 & 0 \\ 
0 & 0 & 0 & 0 & 1 & 0 & -1 & -1 & -1 & 0 & 0 & 1 & 0 & 0 & 1 & 0 \\
0 & 1 & 0 & 0 & 0 & 0 & 0 & 1 & 0 & -1 & -1 & -1 & 0 & 0 & 0 & 1 \\
0 & 0 & 1 & 0 & 0 & 1 & 0 & 0 & 1 & 0 & 0 & 0 & -1 & -1 & -1 & -1 \\
\end{pmatrix} 
\]

The vector $\tilde{b}$, where $\tilde{b}_i$ represents the demand or supply at node $i$, can be computed by recording the given demand or supply at each node and is given by $\tilde{b}=(9,-13,15,-11,0)$. As discussed in the description of the mass balance constraint in Section \ref{subsec:Problem-statement-for-networks}, we compute the full row rank node-arc incidence matrix $A$ by removing the 5th row of $\tilde{A}$
and the resource vector $b$ by removing the last component of $\tilde{b}$, \emph{i.e.}, $b=(9,-13,15,-11)$. So $A$ has 4 rows and 16 columns, and it is:

\[ 
A=   
\begin{pmatrix}
-1 & -1 & -1 & 1 & 0 & 0 & 1 & 0 & 0 & 1 & 0 & 0 & 1 & 0 & 0 & 0 \\ 
1 & 0 & 0 & -1 & -1 & -1 & 0 & 0 & 0 & 0 & 1 & 0 & 0 & 1 & 0 & 0 \\ 
0 & 0 & 0 & 0 & 1 & 0 & -1 & -1 & -1 & 0 & 0 & 1 & 0 & 0 & 1 & 0 \\
0 & 1 & 0 & 0 & 0 & 0 & 0 & 1 & 0 & -1 & -1 & -1 & 0 & 0 & 0 & 1 \\
\end{pmatrix},
\]
where we associate player $i$ with the $i$th column of
$A$. 

\paragraph*{Reindexing the variables.} $\;$To obtain the notational advantage as discussed in the comment after Lemma \ref{28Dec2018}, we next rearrange the columns of
the matrix $A$ so that the unimodular basis matrix comprising of the last 4 columns of $A$ constitutes the
first $4$ columns in the new $A$ \emph{i.e.}, if $A=[A_{1}\mid A_{2}\mid\ldots\mid A_{4}\mid A_{5}\mid\ldots\mid A_{16}]$
then $\det([A_{1}\mid A_{2}\mid\ldots\mid A_{4}])=\pm1$. We reindex
the variables as follows. The last 4 players are reindexed as the first four players, and players 1 to 12 are reindexed as players 5 to 16.  Due to this reindexing procedure, we need to reindex the upper bound vector $u$, where the last 4 indices become the first 4 indices, and the indices 1 to 12 become indices 5 to 16, \emph{i.e.},
\[
u=(5,6,6,10,10,7,11,13,16,12,4,5,6,14,13,15).
\]

Note that indices of $b$ would not change due to the reindexing procedure. We perform our computation on these reindexed variables and then revert the resultant efficient Pareto optimal solutions to the original indexing.

\paragraph*{Problem transformation.}$\;$ Now we are in a position to transform the problem by following the methodology described in Section \ref{sec:transforming_the_game}. First, we decouple the optimization problems for the last 12 players. In order to achieve that, we follow the constructive proofs of Lemma \ref{thm:unimod_mat_ext} and Theorem \ref{thm:transform-x-to-z} to arrive at the following representation of the variable $x$ in terms of the new variable $z$:
%We have from Theorem \ref{thm:transform-x-to-z}
%\[
%U=\begin{bmatrix}I_{4} & -A_{[1:4,5:16]}\\
%0_{12\times4} & I_{12}.
%\end{bmatrix}.
%\]
\begin{equation}
x=\left(\begin{array}{c}
x_{1}\\
x_{2}\\
x_{3}\\
x_{4}\\
x_{5}\\
x_{6}\\
x_{7}\\
x_{8}\\
x_{9}\\
x_{10}\\
x_{11}\\
x_{12}\\
x_{13}\\
x_{14}\\
x_{15}\\
x_{16}
\end{array}\right)=\left(\begin{array}{c}
z_{1}+z_{2}+z_{3}-z_{4}-z_{7}-z_{10}+9\\
-z_{1}+z_{4}+z_{5}+z_{6}-z_{11}-13\\
-z_{5}+z_{7}+z_{8}+z_{9}-z_{12}+15\\
-z_{2}-z_{8}+z_{10}+z_{11}+z_{12}-11\\
z_{1}\\
z_{2}\\
z_{3}\\
z_{4}\\
z_{5}\\
z_{6}\\
z_{7}\\
z_{8}\\
z_{9}\\
z_{10}\\
z_{11}\\
z_{12}
\end{array}\right),
\label{amar_je_shob_dite_hobe}
\end{equation}
where the advantage of this transformation is that, for players 5 to 16, we
have decoupled univariate optimization problems of the form \eqref{eq:decoupled-opt} in $z$, and by solving these decoupled problems, we can reduce the constraint set significantly.

\paragraph*{Solving the decoupled optimization problems.}$\;$Next we solve the aforementioned decoupled univariate optimization problems for
the last 12 players by following the procedure described in the last paragraph of Subsection \ref{decouple_28th_Dec}. The solution set is given by Table \ref{tab:D-for-network}; in this table, the sets of optimal solutions for players $5,6,\ldots,16$ are given by $D_1,D_2,\ldots,D_{12}$, respectively. Solving these 12 decoupled optimization problems immediately reduces the constraint set into a much smaller set, which we define as
$D=\bigtimes_{i=1}^{12}D_{i}$.
\begin{table}
\begin{tabular}{|c|c|}
\hline 
$D_{1}$ & $\{1\}$\tabularnewline
\hline 
$D_{2}$ & $\{3\}$\tabularnewline
\hline 
$D_{3}$ & $\{5\}$\tabularnewline
\hline 
$D_{4}$ & $\{4,6\}$\tabularnewline
\hline 
$D_{5}$ & $\{7,11\}$\tabularnewline
\hline 
$D_{6}$ & $\{10\}$\tabularnewline
\hline 
$D_{7}$ & $\{2\}$\tabularnewline
\hline 
$D_{8}$ & \{1\}\tabularnewline
\hline 
$D_{9}$ & $\{3\}$\tabularnewline
\hline 
$D_{10}$ & $\{7\}$\tabularnewline
\hline 
$D_{11}$ & $\{7\}$\tabularnewline
\hline 
$D_{12}$ & $\{4,5,6,10,11\}$\tabularnewline
\hline 
\end{tabular}\caption{$D$ for the network under consideration \label{tab:D-for-network}}
\end{table}

\begin{table}
\begin{tabular}{|c|c|c|c|c|c|c|c|c|c|c|c|c|}
\hline 
Elements & $z_{1}$ & $z_{2}$ & $z_{3}$ & $z_{4}$ & $z_{5}$ & $z_{6}$ & $z_{7}$ & $z_{8}$ & $z_{9}$ & $z_{10}$ & $z_{11}$ & $z_{12}$\tabularnewline
\hline 
\hline 
1 & 1 & 3 & 5 & 4 & 11 & 10 & 2 & 1 & 3 & 7 & 7 & 4\tabularnewline
\hline 
2 & 1 & 3 & 5 & 4 & 11 & 10 & 2 & 1 & 3 & 7 & 7 & 5\tabularnewline
\hline 
3 & 1 & 3 & 5 & 4 & 11 & 10 & 2 & 1 & 3 & 7 & 7 & 6\tabularnewline
\hline 
4 & 1 & 3 & 5 & 6 & 11 & 10 & 2 & 1 & 3 & 7 & 7 & 4\tabularnewline
\hline 
5 & 1 & 3 & 5 & 6 & 11 & 10 & 2 & 1 & 3 & 7 & 7 & 5\tabularnewline
\hline 
6 & 1 & 3 & 5 & 6 & 11 & 10 & 2 & 1 & 3 & 7 & 7 & 6\tabularnewline
\hline 
\end{tabular}\caption{\label{tab:Table-listing-values-F}Table listing elements of $\mathcal{F}$}
\end{table}
\paragraph*{Consensus reformulation for the first 4 players.}$\;$ Now we are in a position to  transform the optimization problems for
the first 4 players using \eqref{eq:x-in-z} and consensus constraints as described in Subsection \ref{sec:Transforming-the-m-players}. By following the straightforward procedure of Subsection \ref{sec:Transforming-the-m-players}, we arrive at 4 optimization problems of the form \eqref{eq:final-optimization-model-3}, where each of the players 1 to 4 is optimizing its cost function over a \emph{common constraint set} denoted by $\mathcal{F}$. Next we compute the points in $\mathcal{F}$.

\paragraph*{Computing the points in $\mathcal{F}.$}$\;$ To compute the points in $\mathcal{F}$, we follow the methodology described in Subsection \ref{sec:Nonemptyness-of-}. Finding the points in  $\mathcal{F}$ requires computing the reduced Groebner basis with respect to lexicographic ordering, denoted by $G_{\textrm{reduced},\succ_{\textrm{lex}}}$. We compute $G_{\textrm{reduced},\succ_{\textrm{lex}}}$ using the \texttt{GroebnerBasis} function in \texttt{Wolfram Mathematica 10}, and we find it to be not equal to $\{1\}$. Hence $\mathcal{F}$ is nonempty due to Lemma \ref{choke_amar_trishna}. Next we compute the points in $\mathcal{F}$ by following Algorithm \ref{alg:extracting_mathcal_F-2}; the list of computed points is given by Table \ref{tab:Table-listing-values-F}.

\paragraph*{Computing the efficient Pareto optimal points.}$\;$The final step is  computing the efficient
Pareto optimal points  from $\mathcal{F}$ by executing  the three steps described in Algorithm \ref{alg:other-m-players-1}. After applying Algorithm \ref{alg:other-m-players-1}, we arrive at two efficient Pareto optimal points in variable $z$ as follows:
\[
(1,3,5,4,11,10,2,1,3,7,7,5),
\]
 and 
\[
(1,3,5,6,11,10,2,1,3,7,7,6).
\]

Using the relationship between $x$ and $z$ in \eqref{amar_je_shob_dite_hobe}, we can express these efficient Pareto optimal points in variable $x$ as follows:

\[
(5,4,5,4,1,3,5,4,11,10,2,1,3,7,7,5),
\]
 and 
\[
(3,6,4,5,1,3,5,6,11,10,2,1,3,7,7,6).
\]

Note that these efficient  Pareto optimal points above are associated with the reindexed $x$. By reversing this reindexing, where indices 5 to 16 will be indices 1 to 12, and indices 1 to 4 will be indices 13 to 16, we arrive at the efficient Pareto optimal points in our original variable $x$ as follows:

\[
(1,3,5,4,11,10,2,1,3,7,7,5, 5,4,5,4),
\]
 and 
\[
(1,3,5,6,11,10,2,1,3,7,7,6,3,6,4,5).
\]

\section{Conclusions} \label{sec:Conclusion}

\textcolor{black}{In this paper, we have proposed a multi-player extension
of the minimum cost flow problem inspired by a multi-player transportation
problem. We have associated
one player with each arc of a directed connected network, each trying
to minimize its cost subject to the network flow constraints. The
cost can be any general nonlinear function, and the flow through each
arc is integer-valued. In this multi-player setup, we have presented
algorithms to compute an }\textcolor{black}{\emph{efficient Pareto
optimal point}}\textcolor{black}{, which is a good compromise solution
between the utopian vector optimality and the generic Pareto optimality
and} is a Nash equilibrium\textcolor{black}{{} if our problem
is transformed into an $n$-person finite static game in normal form. }

Some concluding remarks on the limitations of our methodology are
as follows. \emph{First}, at the heart of our methodology is the transformation
provided by Theorem \ref{thm:transform-x-to-z}, which decouples the
optimization problems for the last $n-m$ players. Each of these decoupled
optimization problems is univariate over a discrete interval and
is easy to solve. This can potentially allow us to work in a much
smaller search space. So if we have a system where $n-m>m\Leftrightarrow m<\frac{n}{2}$,
then it will be convenient from a numerical point of view. \emph{Second},
computation of Pareto optimal points depends on determining the points
in $\mathcal{F}$ using Groebner basis. Calculating Groebner basis
can be numerically challenging for large system \citep[pages 111-112]{Cox2007},
though significant speed-up has been achieved  in recent years by computer
algebra packages such as Macaulry2, SINGULAR, FGb, and Mathematica. 

\bibliographystyle{unsrtnat}
%\bibliography{Mnscrpt_pareto_opt}
% contents of the .bbl file starts here

% contents of the .bbl file ends here
%-------------------------------------------

\newpage 

\section*{Appendix}

\textbf{Proof of Lemma 1}

\begin{comment}
\begin{customlemma}{1}\label{eight} 
Assumption \ref{assum:shared_vertex} holds if and only if we can extract a unimodular basis matrix from $A$.
\end{customlemma}
\end{comment}

\noindent Any basic solution of $\mathop{\bf relaxed}P$ can be constructed
as follows. Take $m$ linearly independent columns 
\[
A_{B(1)},A_{B(2)},\ldots,A_{B(m)},
\]
 where $B(1),B(2),\ldots,B(m)$ are the indices of those columns.
Construct the basis matrix 
\[
B=\begin{bmatrix}A_{B(1)} &  & A_{B(2)} &  & \cdots &  & A_{B(m)}\end{bmatrix}.
\]
Set $x_{B}=(x_{B(1)},x_{B(2)},\ldots,x_{B(m)})=B^{-1}b$, which is
called the basic variable in linear programming theory. Set the rest
of the components of $x$ (called the nonbasic variable) to be zero,
\emph{i.e.}, 
\[
x_{NB}=(x_{NB(1)},x_{NB(2)},\ldots,x_{NB(n-m)})=(0,0,\ldots,0).
\]
The resultant $x$ will be a basic solution of $\mathop{\bf relaxed}P$
\citep[Theorem 2.4]{Bertsimas1997}. 

Suppose there is an integer basic solution in $\mathop{\bf relaxed}P$.
Denote that integer basic solution by $\bar{x}$ and the associated
basis matrix by $\bar{B}$. The nonbasic variables are integers (zeros)
in every basic solution, so the basic variables $\bar{x}_{\bar{B}}=\bar{B}^{-1}b$
has to be an integer vector for any integer $b$. Now from Cramer's
rule \citep[Proposition 3.1]{Bertsimas2005}, 
\[
\left(\forall b\in\mathbf{Z}^{m}\right)\quad\left(\bar{x}_{\bar{B}}=\bar{B}^{-1}b\in\mathbf{Z}^{m}\Leftrightarrow\forall i\in[m]\quad\bar{x}_{\bar{B}(i)}=\frac{\det\bar{B}^{i}}{\det\bar{B}}\in\mathbf{Z}\right),
\]
 where $\bar{B}^{i}$ is the same as $\bar{B}$, except the $i$th
column has been replaced with $\bar{b}$. Now $\det\bar{B}^{i}\in\mathbf{Z}$
because $b$ is an integer vector. So having integer basic solution
is equivalent to $\det\bar{B}=\pm1$, \emph{i.e.}, $\bar{B}$ is unimodular.  \hfill    $\blacksquare$ \\

\noindent\textbf{Proof of Lemma 2}

\begin{comment}
\begin{customlemma}{2}
Let, $C=B^{-1}A$ and $d=B^{-1}b$. Then $C\in\mathbf{Z}^{m\times n}$, $d\in\mathbf{Z}^{m}$, and the constraint set $Q$ and $P$ (defined in \eqref{eq:constraint-set-P} and \eqref{eq:equality-constraint-set-Q}, respectively) have the equivalent representation: \begin{align} Q & =\left\{ x\in\mathbf{Z}^{n}\mid Cx=d\right\} ,\\ P & =\left\{ x\in\mathbf{Z}^{n}\mid Cx=d,0\preceq x\preceq u\right\} .
\end{align} 
\end{customlemma}
\end{comment}

\noindent First, note that unimodularity of $B$ is equivalent to unimodularity
of $B^{-1}$, which we show easily as follows. First note that, $\det\left(B^{-1}\right)=\frac{1}{\det B}=\frac{1}{\pm1}=\pm1$.
Each component of $B^{-1}$ is a subdeterminant of $B$ divided by
$\det B$. As $B\in\mathbf{Z}^{m\times m}$ and $\det B=\pm1$, each
component of $B^{-1}$ is an integer. Thus $B^{-1}$ is a unimodular
matrix. Similarly, unimodularity of $B^{-1}$ implies $B$ is unimodular. 

So $C=B^{-1}A\in\mathbf{Z}^{m\times n}$ and $d=B^{-1}b\in\mathbf{Z}^{m}$
as $A$ and $b$ are integer matrices. As multiplying both sides of
a linear system by an invertible matrix does not change the solution,
we have $Ax=b\Leftrightarrow B^{-1}Ax=B^{-1}b\Leftrightarrow Cx=d$.
Thus we arrive at the claim. \hfill    $\blacksquare$ \\

\noindent\textbf{Proof of Lemma 3}

\begin{comment}
\begin{customlemma}{3} \label{lem:The-matrix-C-conversion-1} 
The matrix $C$ as defined in Lemma \ref{lem:Equivalence-of-constraint-set} can be brought into the Hermite normal form $[I|0]$ by elementary integer column operations, more specifically by adding integer multiple of one column to another column. 
\end{customlemma}
\end{comment}

\noindent First, recall that the following operations on a matrix are called
\emph{elementary integer column operations}: (i) adding an integer
multiple of one column to another column, (ii) exchanging two columns,
and (iii) multiplying a column by -1. 

The matrix $C$ is of the form:
\begin{eqnarray*}
C & = & B^{-1}A\\
 & = & B^{-1}\left[B\mid A_{m+1}\mid A_{m+2}\mid\cdots\mid A_{n}\right]\\
 & = & \begin{bmatrix}1 & 0 & \cdots & 0 & C_{1,m+1} & C_{1,m+2} & \cdots & C_{1,n}\\
0 & 1 & \cdots & 0 & C_{2,m+1} & C_{2,m+2} & \cdots & C_{2,n}\\
\vdots & \vdots & \vdots & \vdots & \vdots & \vdots & \vdots & \vdots\\
0 & 0 & \cdots & 1 & C_{m,m+1} & C_{m,m+2} & \cdots & C_{m,n}
\end{bmatrix}
\end{eqnarray*}
Consider the first column of $C$. For all $j=m+1,m+2,\ldots,n$,
we multiply the first column of $C$, $C_{1}=e_{1}$ of $C$ by $-C_{1,j}$
and then add it to $C_{i}$. Thus $C$ is transformed to
\[
\begin{bmatrix}1 & 0 & \cdots & 0 & 0 & 0 & \cdots & 0\\
0 & 1 & \cdots & 0 & C_{2,m+1} & C_{2,m+2} & \cdots & C_{2,n}\\
\vdots & \vdots & \vdots & \vdots & \vdots & \vdots & \vdots & \vdots\\
0 & 0 & \cdots & 1 & C_{m,m+1} & C_{m,m+2} & \cdots & C_{m,n}
\end{bmatrix}
\]
Similarly for column indices, $i=2,3,\ldots,m$, respectively we do
the following. For $j=m+1,m+2,\cdots,n$, we multiply the $i$th column
$e_{i}$ with $-C_{i,j}$ and add it to $C_{j}$. In the end, the Hermite
normal form becomes: 
\[
\begin{bmatrix}1 & 0 & \cdots & 0 & 0 & 0 & \cdots & 0\\
0 & 1 & \cdots & 0 & 0 & 0 & \cdots & 0\\
\vdots & \vdots & \vdots & \vdots & \vdots & \vdots & \vdots & \vdots\\
0 & 0 & \cdots & 1 & 0 & 0 & \cdots & 0
\end{bmatrix}.
\]

The steps describing the process in Lemma \ref{lem:The-matrix-C-conversion}
can be summarized by Algorithm \ref{alg:mat_to_Hermite_normal_form},
which we are going to use to prove Lemma 4.  \hfill    $\blacksquare$ \\

\begin{algorithm} 
 \caption{Converting $C$ to $[I|0]$}
\label{alg:mat_to_Hermite_normal_form}
\begin{algorithmic}[1] 
 \Procedure {Converting $C$ to $[I|0]$}{}
\For{$i:= 1,2,\ldots, m$} 
\For{$j:= m+1,m+2,\ldots,n$} 
\State $C_j := C_j-C_{i,j} e_i$ 
\EndFor 
\EndFor 
 \EndProcedure
\end{algorithmic} 
\end{algorithm}

\noindent\textbf{Proof of Lemma 4}

\begin{comment}
\begin{customlemma}{4} \label{thm:unimod_mat_ext} 
There exists a unimodular matrix $U$ such that $CU=\left[I\mid0\right]$. 
\end{customlemma}
\end{comment}

\noindent If we multiply column $C_{i}$ of a matrix $C$ with an integer factor
$\gamma$ and subsequently add it to another column $C_{j}$, then it is equivalent
to right multiplying the matrix $C$ with a matrix $I+\gamma I_{ij}$
(recall that the matrix $I_{ij}$ has a one in $(i,j)$th position
and zero everywhere else). Note that $I+\gamma I_{ij}$ is a triangular
matrix with diagonal entries being one, $\gamma$ being on the $(i,j)$th
position, and zero everywhere else. As the determinant of a triangular
matrix is the product of its diagonal entries, $\det(I+\gamma I_{ij})=1$.
So $I+\gamma I_{ij}$ is a unimodular matrix. Furthermore, step 4 of the procedure
above to convert $C$ to Hermite normal form, \emph{i.e.}, $C_{j}=C_{j}-C_{i,j}e_{i}$
is equivalent to left multiplying the current matrix with $I-C_{i,j}I_{ij}$.
So the inner loop of the procedure above over $j=m+1,m+2,\ldots,n$ (lines
2-4) can be achieved by left multiplying the current matrix with 
\[
U_{i}=\prod_{j=m+1}^{n}(I-C_{ij}I_{ij})=(I-C_{i,m+1}I_{i,m+1})(I-C_{i,m+2}I_{i,m+2})\cdots(I-C_{i,n}I_{i,n})
\]
As each of the matrices in the product is a unimodular matrix and
determinant of multiplication of square matrices of same size is equal
to multiplication of determinant of those matrices, we have $\det(U_{i})=1$.
So $U_{i}$ is a unimodular matrix. Structurally the $i$th row of
$U_{i}$, denoted by $u_{i}^{T}$, has a $1$ on $i$th position,
has $-C_{i,j}$ on $j$th position for $j=m+1,m+2,\ldots,n$, and
zero everywhere else. Any other $k$th row $(k\neq i)$ of $U_{i}$
is $e_{k}^{T}$. So we can convert $C$ to its Hermite normal form
by repeatedly left multiplying $C$ with $U_{1},U_{2},\ldots,U_{m}$, respectively. This is equivalent to left multiplying $C$ with one
single matrix $U=\prod_{i=1}^{m}U_{i}$. The final matrix $U$ is
unimodular as it is multiplication of unimodular matrices.  \hfill    $\blacksquare$ \\

\noindent\textbf{Proof of Theorem 1}

\begin{comment}
\begin{customthm} {1} \label{thm:transform-x-to-z}
The constraint set $Q$ defined in \eqref{eq:constraint_set_v2} is nonempty and for any vector $x$ can be maximally decomposed in terms of a new variable $z$ as \begin{equation*} x\in Q\Leftrightarrow\exists z\in\mathbf{Z}^{n-m}\quad x=\begin{bmatrix}d_{1}-h_{1}^{T}z\\ d_{2}-h_{2}^{T}z\\ \vdots\\ d_{m}-h_{m}^{T}z\\ z_{1}\\ \vdots\\ z_{n-m} \end{bmatrix}, \end{equation*}  where $d_{i}$ is the $i$th component of $d=B^{-1}b$, and $h_{i}^{T}\in\mathbf{Z}^{n-m}$ is the $i$th row of $B^{-1}A_{[1:m,m+1:n]}$.  
\end{customthm}
\end{comment}

\noindent Let $y=U^{-1}x$. As $U$ is unimodular, $U^{-1}$ is also unimodular. Hence $x\in\mathbf{Z}^{n}\Leftrightarrow y\in\mathbf{Z}^{n}$.
Let $y=(y_{1},y_{2})$, where $y_{1}\in\mathbf{Z}^{m}$ and $y_{2}\in\mathbf{Z}^{n-m}$.
Then

\begin{align*}
 & Q\neq\emptyset\\
\Leftrightarrow & \exists x\in\mathbf{Z}^{n}\;\left(Cx=d\right)\\
\Leftrightarrow & \exists y\in\mathbf{Z}^{n}\;\left(CUy=\begin{bmatrix}I & 0\end{bmatrix}\begin{bmatrix}y_{1}\\
y_{2}
\end{bmatrix}=y_{1}=d\right).
\end{align*}

As $d=B^{-1}b\in\mathbf{Z}^{m}$ (Lemma \ref{lem:Equivalence-of-constraint-set}),
by taking $y=\begin{bmatrix}y_{1}\\
y_{2}
\end{bmatrix}=\begin{bmatrix}B^{-1}b\\
z
\end{bmatrix}\in\mathbf{Z}^{n}$, where $z\in\mathbf{Z}^{n-m}$, we can satisfy the condition above.
Thus $Q$ is nonempty. 

Now for any $x$, we can maximally decompose it in terms of $z$ as:
\begin{align*}
x\in Q\Leftrightarrow x=Uy & =\begin{bmatrix}I_{m} & -B^{-1}A_{[1:m,m+1:n]}\\
0_{n-m\times m} & I_{n-m}.
\end{bmatrix}\begin{bmatrix}B^{-1}b\\
z
\end{bmatrix}\\
 & =\begin{bmatrix}B^{-1}\left(b-A_{[1:m,m+1:n]}z\right)\\
z
\end{bmatrix}\\
 & =\begin{bmatrix}d_{1}-h_{1}^{T}z\\
d_{2}-h_{2}^{T}z\\
\vdots\\
d_{m}-h_{m}^{T}z\\
z_{1}\\
\vdots\\
z_{n-m}
\end{bmatrix},
\end{align*}
 where the last $n-m$ variables are completely decoupled in terms
of $z$. Note that the decomposition is maximal by construction and
cannot be extended any further in general cases. \hfill    $\blacksquare$ \\

\noindent\textbf{Proof of Theorem 3}

\begin{comment}
\begin{customthm}{3}
The set $\mathcal{F}$ is nonempty if and only  \[ G_{\textrm{reduced},\succ}\ne\{1\}, \]  where $G_{\textrm{reduced},\succ}$ is the reduced Groebner basis of $\mathbf{ideal}\left\{ q_{1},\ldots,q_{m},r_{1},\ldots,r_{n-m}\right\} $ with respect to any ordering. 
\end{customthm}
\end{comment}

\noindent The proof sketch is as follows. The elements of $\mathcal{F}$ are
the solution of the polynomial system: 
\begin{align}
 & \left(\forall i\in[m]\right)\quad q_{i}(z)=0\nonumber \\
 & \left(\forall j\in[n-m]\right)\quad r_{j}(z)=0.\label{eq:poly_sys_2-2}
\end{align}
We prove in two steps. In step $1$, we show that the polynomial system
\eqref{eq:poly_sys_2-2} is feasible if and only if 
\[
1\notin\mathbf{ideal}\left\{ q_{1},\ldots,q_{m},r_{1},\ldots,r_{n-m}\right\} .
\]
Then in step $2$, we show that, $1\notin\mathbf{ideal}\left\{ q_{1},\ldots,q_{m},r_{1},\ldots,r_{n-m}\right\} $
is equivalent to $G_{\textrm{reduced},\succ}\ne\{1\}$. \paragraph{\textcolor{blue}{Step 1. Polynomial system \eqref{eq:poly_sys_2-2} is feasible if and only if $1\notin\mathbf{ideal}\left\{ q_{1},\ldots,q_{m},r_{1},\ldots,r_{n-m}\right\}.$}}
 We prove necessity first and then sufficiency. \newline ($\Rightarrow$)

Assume system \eqref{eq:poly_sys_2-2} is feasible, but $1\in\mathbf{ideal}\left\{ q_{1},\ldots,q_{m},r_{1},\ldots,r_{n-m}\right\} $.
That means 

\begin{equation}
\left(\exists h_{1},\ldots,h_{m},s_{1},\ldots,s_{n-m}\in\mathbf{C}[z]\right)\left(\forall z\in\mathbf{C}^{n-m}\right)\quad1=\sum_{i=1}^{m}h_{i}(z)q_{i}(z)+\sum_{i=1}^{n-m}s_{i}(z)r_{i}(z),\label{eq:Farka_type_lemma_1}
\end{equation}
 and 

\begin{align}
\left(\exists\bar{z}\in\mathbf{Z}{}^{n-m}\right)\bigg( & \left(\forall i\in[m]\right)\quad q_{i}(\bar{z})=0,\nonumber \\
 & \left(\forall i\in[n-m]\right)\quad r_{i}(\bar{z})=0\bigg).\label{eq:Farka_type_lemma_2}
\end{align}
Putting $z=\bar{z}$ in \eqref{eq:Farka_type_lemma_1} and
using \eqref{eq:Farka_type_lemma_2} we get $1=0$, so we
have a contradiction. \newline ($\Leftarrow$)

We want to show that if $1\not\in\mathbf{ideal}\{q_{1},\cdots,q_{m},r_{1},\cdots,r_{n-m}\}$,
then the polynomial system \eqref{eq:poly_sys_2-2} is feasible. We
prove the contrapositive again: if the polynomial system \eqref{eq:poly_sys_2-2}
is infeasible in integers, then $1\in\mathbf{ideal}\{q_{1},\cdots,q_{m},r_{1},\cdots,r_{n-m}\}$. 

First, we show that, \emph{feasibility of the system in $\mathbf{C}^{n-m}$
is equivalent to feasibility in $\mathbf{Z}^{n-m}$. }As $\mathbf{Z}^{n-m}\subset\mathbf{C}^{n-m}$,
if the polynomial system is infeasible in $\mathbf{C}^{n-m}$, it
is infeasible in $\mathbf{Z}^{n-m}$. Also, if the polynomial system
is infeasible in $\mathbf{Z}^{n-m}$, then it will be infeasible in
$\mathbf{C}^{n-m}$. We show this by contradiction. Assume the system
system is infeasible in $\mathbf{Z}^{n-m}$ but feasible in $\mathbf{C}^{n-m}$
\emph{i.e.}, $\mathbf{C}^{n-m}\setminus\mathbf{Z}^{n-m}$. Let that
feasible solution be $\tilde{z}\in\mathbf{C}^{n-m}\setminus\mathbf{Z}^{n-m}$,
so there is at least one component of it (say $\tilde{i}$) such that
$\tilde{z}_{\tilde{i}}\in\mathbf{C}\setminus\mathbf{Z}$. Now
\[
r_{\tilde{i}}(\tilde{z})=(\tilde{z}_{\tilde{i}}-z_{i,1})(\tilde{z}_{\tilde{i}}-z_{i,2})\ldots(\tilde{z}_{\tilde{i}}-z_{i,p_{i}}),
\]
where, by construction, each of the elements of the set $D_{i}=\left\{ z_{i,1},z_{i,2},\ldots,z_{i,p_{i}}\right\} $
are integers and different from each other, so each component in the
product $r_{\tilde{i}}(\tilde{z})$ are nonzero complex numbers with
the absence of complex conjugates. So $r_{\tilde{i}}(\tilde{z})\neq0$,
which is a contradiction. 

If the polynomial system is infeasible in $\mathbf{C}^{n-m}$, then
it is equivalent to saying that, 
\[
V(q_{1},\ldots,q_{m},r_{1},\ldots,r_{n-m})=\emptyset,
\]
where $V$ has been defined in \eqref{eq:aff_variety_def}.
Then using the \emph{Weak Nullstellensatz}, we have
\[
\mathbf{ideal}\{q_{1},\ldots,q_{m},r_{1},\ldots,r_{n-m}\}=\mathbf{C}[z_{1},z_{2},\ldots,z_{n-m}].
\]
As $1\in\mathbf{C}[z_{1},z_{2},\ldots,z_{n-m}]$, this means $1\in\mathbf{ideal}\{q_{1},\ldots,q_{m},r_{1},\ldots,r_{n-m}\}$.
\paragraph{\textcolor{blue}{Step 2. Now we show $1\notin\mathbf{ideal}\left\{ q_{1},\ldots,q_{m},r_{1},\ldots,r_{n-m}\right\} $ is equivalent to $G_{\textrm{reduced},\succ}\ne\{1\}.$}} 

We have shown that feasibility of the system in $\mathbf{C}^{n-m}$
is equivalent to feasibility in $\mathbf{Z}^{n-m}$. As a result,
we can work over complex numbers, which is a algebraically closed
field. This allows us to apply \emph{consistency algorithm} \citep[page 172]{Cox2007},
which states that $1\notin\mathbf{ideal}\left\{ q_{1},\ldots,q_{m},r_{1},\ldots,r_{n-m}\right\} $
if and only if $G_{\textrm{reduced},\succ}\ne\{1\}$ .  \hfill    $\blacksquare$ \\

\noindent\textbf{Proof of Lemma 5}

\begin{comment}
\begin{customlemma}{5}  \label{lem:groebner-basis-lemma}Suppose $G_{\textrm{reduced},\succ}\ne\{1\}$. Then $\mathcal{F}=V\left(G_{\textrm{reduced},\succ_{\textrm{lex}}}\right)\neq\emptyset.$ \end{customlemma}
\end{comment}

\noindent By Theorem \ref{lem:poly_sys_ideal_lemma}, we have $\mathcal{F}\neq\emptyset$.
So from the definition of affine variety in \eqref{eq:aff_variety_def}
and \eqref{eq:mathcal-F} we can write $\mathcal{F}$ as: 
\begin{align*}
 & \mathcal{F}=V\left(q_{1},\ldots,q_{m},r_{1},\ldots,r_{n-m}\right)\cap\mathbf{Z}^{n-m}.
\end{align*}
From the equation above and \eqref{eq:variety equality} we have,
\begin{eqnarray*}
\mathcal{F} & = & V\left(q_{1},\ldots,q_{m},r_{1},\ldots,r_{n-m}\right)\\
 & \neq & \emptyset.
\end{eqnarray*}
 By definition of basis, we have $\mathbf{ideal}\left\{ q_{1},\ldots,q_{m},r_{1},\ldots,r_{n-m}\right\} =\mathbf{ideal}\left\{ G_{\textrm{reduced},\succ_{\textrm{lex}}}\right\} $,
which implies 
\[
V\left(q_{1},\ldots,q_{m},r_{1},\ldots,r_{n-m}\right)=V\left(G_{\textrm{reduced},\succ_{\textrm{lex}}}\right)
\]
 due to \citep[Proposition 2,][page 32]{Cox2007}. So $\mathcal{F}=V\left(G_{\textrm{reduced},\succ_{\textrm{lex}}}\right).$ \hfill    $\blacksquare$ \\

\noindent\textbf{Proof of Lemma 6}

\begin{comment}
\begin{customlemma}{6}  Algorithm \ref{alg:extracting_mathcal_F-2} correctly calculates all the points in $\mathcal{F}$, when it is nonempty. \label{lem:Algorithm--correctly} \end{customlemma}
\end{comment}

\noindent Using Theorem \ref{lem:poly_sys_ideal_lemma}, we have $G_{\textrm{reduced},\succ_{\textrm{lex}}}\neq\emptyset$.
So by the elimination theorem \citep[Theorem 2,][page 116]{Cox2007}
$V\left(G_{n-m-1}\right)$ is nonempty and will contain the list of all possible
$z_{n-m}$ coordinates for the points in $\mathcal{F}$. As $G_{\textrm{reduced},\succ_{\textrm{lex}}}\neq\emptyset$,
when moving from one step to the next, not all the affine varieties
associated with the univariate polynomials (after replacing the previous
coordinates into the elimination ideal) can be empty due to the extension
theorem \citep[Theorem 3,][page 118]{Cox2007}. Using this logic repeatedly,
the final step will give us $\mathcal{F}=V(G_{\textrm{reduced},\succ_{\textrm{lex}}})\neq\emptyset$.    \hfill $\blacksquare$ \\

\noindent\textbf{Proof of Lemma 7}

\begin{comment}
\begin{customlemma}{7} \label{lem:subseteq} In Algorithm \ref{alg:other-m-players-1}, for all $i\in[m-1]$ we have $\mathcal{F}_{s_{i+1}}^{*}\subseteq\mathcal{F}_{s_{i}}^{*}\subseteq\mathcal{F}$. \end{customlemma}
\end{comment}

\noindent Follows from \eqref{eq:mcFsi}, \eqref{eq:xi-inv} and \eqref{eq:Xsi1}.        \hfill $\blacksquare$ \\

\noindent\textbf{Proof of Lemma 8}

\begin{comment}
\begin{customlemma}{8}\label{lem:optimal-x-z}In Algorithm \ref{alg:other-m-players-1} for any $i\in[m]$, $x_{s_{i}}\in X_{s_{i}}^{*}$ if and only if $z^{*}\in\mathcal{F}_{s_{i}}^{*}$ . Furthermore, $z^{*}\in\mathcal{F}_{s_{i}}^{*}$ solves the following optimization problem  \begin{align*}  & \begin{aligned} & \mathrm{minimize}_{z} &  & f_{s_{i}}(d_{s_{i}}-h_{s_{i}}^{T}z)\\  & \mathrm{subject\;to} &  & z\in\mathcal{F}_{s_{i-1}}^{*}, \end{aligned} \end{align*} for all $i\in[2:m]$. \end{customlemma}
\end{comment}

\noindent For any $i\in[m]$, 
\begin{align*}
 & x_{s_{i}}\in X_{s_{i}}^{*}\\
\Leftrightarrow & d_{s_{i}}-h_{s_{i}}^{T}z\in X_{s_{i}}^{*}\\
\Leftrightarrow & z\in\mathcal{F}_{s_{i}}^{*},
\end{align*}
 where the second line follows from \eqref{eq:x-in-z}. So
\begin{align*}
 & \left(\begin{aligned} & \mathrm{minimize}_{x_{s_{i}}} &  & f_{s_{i}}(x_{s_{i}})\\
 & \mathrm{subject\;to} &  & x_{s_{i}}\in X_{s_{i}}
\end{aligned}
\right)\\
= & \left(\begin{aligned} & \mathrm{minimize}_{x_{s_{i}}} &  & f_{s_{i}}(d_{s_{i}}-h_{s_{i}}^{T}z)\\
 & \mathrm{subject\;to} &  & z\in\mathcal{F}_{s_{i-1}}^{*}
\end{aligned}
\right),
\end{align*}
 where the second line follows from \eqref{eq:Xsi1} in Algorithm
\ref{alg:other-m-players-1}.   \hfill $\blacksquare$ \\

\noindent\textbf{Proof of Lemma 9}

\begin{comment}
\begin{customlemma}{9}
Suppose $\mathcal{F}\neq\emptyset$. Then in Algorithm \ref{alg:other-m-players-1}, $\mathcal{F}_{s_{i}}^{*}$ is nonempty for any $i\in[m]$.
\end{customlemma}
\end{comment}

\noindent As $\mathcal{F}\neq\emptyset$, $X_{s_{1}}:=\left\{ d_{s_{1}}-h_{s_{1}}^{T}z^{(s_{1})}\mid z^{(s_{1})}\in\mathcal{F}\right\} \neq\emptyset$.
Assume, for $i\in[m]$, we have $\mathcal{F}_{s_{i}}^{*}\neq\emptyset$.
Then $X_{s_{i+1}}:=\left\{ d_{s_{i+1}}-h_{s_{i+1}}^{T}z\mid z\in\mbox{\ensuremath{\mathcal{F}}}_{s_{i}}^{*}\right\} =d_{s_{i+1}}-h_{s_{i+1}}^{T}\mathcal{F}_{s_{i}}^{*}\neq\emptyset$.
The subsequent optimization problem is 
\[
\begin{aligned} & \mathrm{minimize}_{x_{s_{i+1}}} &  & f_{s_{i+1}}(x_{s_{i+1}})\\
 & \mathrm{subject\;to} &  & x_{s_{i+1}}\in X_{s_{i+1}}.
\end{aligned}
\]
As we are optimizing over a finite and countable set, a minimizer
will exist. So $X_{s_{i+1}}^{*}\neq\emptyset$. Hence $\mathcal{F}_{s_{i}}^{*}=\bigcup_{x_{i}\in X_{s_{i}}^{*}}(X_{s_{i}}^{*})^{-1}(x_{i})\neq\emptyset$.
So for any $i=1,\ldots,m$, we have $\mathcal{F}_{s_{i}}^{*}$ nonempty. \hfill $\blacksquare$ \\

\noindent\textbf{Proof of Theorem 4}

\begin{comment}
\begin{customthm}{4}
For any $z^{*}\in\mathcal{F}_{s_{m}}^{*}$,  \begin{equation} x^{*}=(d_{1}-h_{1}^{T}z^{*},\ldots,d_{m}-h_{m}^{T}z^{*},z_{1}^{*},\ldots,z_{n-m}^{*})\label{eq:x*z*} \end{equation}  is a Pareto optimal point.
\end{customthm}
\end{comment}

\noindent We want to show that $x^{*}$ is feasible, and for every $x\in P, i\in[n]$ , if we have $f_{i}(x_{i}^{*})\geq f_{i}(x_{i})$, then
for every $j\in[n]$, $f_{j}(x_{j}^{*})=f_{j}(x_{j})$. Using \eqref{eq:x*z*} and \eqref{eq:x-in-z}, we can translate the Pareto
optimality condition in $z$ as follows. Consider a $z\in\mathbf{Z}^{n-m}$
such that 
\begin{eqnarray}
(0,\ldots,0)\preceq\left(\left(d_{i}-h_{i}^{T}z\right)_{i=1}^{m}\right),z)\preceq & (u_{1},\ldots,u_{n}),\label{eq:given-Pareto-proof}
\end{eqnarray}
and suppose 

\begin{align}
\left(\left(f_{i}(d_{i}-h_{i}^{T}z^{*})_{i=1}^{m}\right),\left(f_{m+i}(z_{i}^{*})\right)_{i=1}^{n-m}\right)\succeq & \left(\left(f_{i}(d_{i}-h_{i}^{T}z)_{i=1}^{m}\right),\left(f_{m+i}(z_{i})\right)_{i=1}^{n-m}\right).\label{eq:antecedent-pareto-proof}
\end{align}
 Then we want to show that: 
\begin{align}
\left(\left(f_{i}(d_{i}-h_{i}^{T}z^{*})_{i=1}^{m}\right),\left(f_{m+i}(z_{i}^{*})\right)_{i=1}^{n-m}\right)= & \left(\left(f_{i}(d_{i}-h_{i}^{T}z)_{i=1}^{m}\right),\left(f_{m+i}(z_{i})\right)_{i=1}^{n-m}\right)\label{eq:goal-pareto-z}
\end{align}

Let us start with the last $n-m$ rows. As, $z^{*}\in\mathcal{F}_{s_{m}}^{*}\subseteq\mathcal{F}\subseteq D$
and by construction, $D=\bigtimes_{i=1}^{n-m}D_{i}$, where any element
of $D_{i}$ is a minimizer of \eqref{eq:decoupled-opt}, so 
\[
\left(f_{m+i}(z_{i}^{*})\right)_{i=1}^{n-m}\succeq\left(f_{m+i}(z_{i})\right)_{i=1}^{n-m}
\]
implies 
\[
\left(f_{m+i}(z_{i}^{*})\right)_{i=1}^{n-m}=\left(f_{m+i}(z_{i})\right)_{i=1}^{n-m}.
\]

In the subsequent steps, it suffices to confine $z\in D$, as otherwise
last $n-m$ inequalities of \eqref{eq:antecedent-pareto-proof}
will be violated. Now let us consider the first $m$ inequalities
of \eqref{eq:antecedent-pareto-proof} As discussed in Section
\ref{sec:Transforming-the-m-players}, any $z$ is in $D$ which obeys
the inequalities $0\leq d_{i}-h_{i}^{T}z\leq u_{i}$ for $i=1,\ldots,m$;
this is equivalent to $z\in\mathcal{F}\subseteq D$. Consider, $s_{1}\in[m]$
. Lemmas \ref{lem:subseteq} and \ref{lem:optimal-x-z} implies that
$z^{*}$ solves the following optimization problem $\min_{z}\{f_{s_{1}}(d_{s_{1}}-h_{s_{1}}^{T}z)\mid z\in\mathcal{F}\}=\min_{x_{s_{1}}}\{f_{s_{1}}(x_{s_{1}})\mid x_{s_{1}}\in X_{s_{1}}\}$,
which has solution $x_{s_{1}}^{*}\in X_{s_{1}}^{*}\Leftrightarrow z^{*}\in\mathcal{F}_{s_{1}}^{*}\supseteq\mathcal{F}_{s_{m}}^{*}$.
So $f_{s_{1}}(d_{s_{1}}-h_{s_{1}}^{T}z)\leq f_{s_{1}}(d_{s_{1}}-h_{s_{1}}^{T}z^{*})$
implies $f_{s_{1}}(d_{s_{1}}-h_{s_{1}}^{T}z)=f_{s_{1}}(d_{s_{1}}-h_{s_{1}}^{T}z^{*})$
and $z\in\mathcal{F}_{s_{1}}^{*}$.

Consider $s_{2}\in[m]\setminus\{s_{1}\}$. First, note that $z\in\mathcal{F}_{s_{1}}^{*}$,
otherwise $f_{s_{1}}(d_{s_{1}}-h_{s_{1}}^{T}z)\leq f_{s_{1}}(d_{s_{1}}-h_{s_{1}}^{T}z^{*})$
will not hold. Now the $x_{s_{2}}^{*}$ associated with $z^{*}$ solves
the following optimization problem $\min_{x_{s_{2}}}\{f_{s_{2}}(x_{s_{2}})\mid x_{s_{2}}\in X_{s_{2}}\}=\min_{z}\{f_{s_{2}}(d_{s_{2}}-h_{s_{2}}^{T}z)\mid z\in\mathcal{F}_{s_{1}}^{*}\}$,
where an optimal solution to the first line will be in $X_{s_{2}}^{*}$
and the optimal solution to the second line will be in $\mathcal{F}_{s_{2}}^{*}$
(Lemma \ref{lem:optimal-x-z}). So combining $f_{s_{2}}(d_{s_{2}}-h_{s_{2}}^{T}z)\leq f_{s_{2}}(d_{s_{2}}-h_{s_{2}}^{T}z^{*})$
and $z\in\mathcal{F}_{s_{1}}^{*}$ implies $f_{s_{2}}(d_{s_{2}}-h_{s_{2}}^{T}z)=f_{s_{2}}(d_{s_{2}}-h_{s_{2}}^{T}z^{*}).$
Repeating a similar argument for $i=s_{3},s_{4},\ldots,s_{m}$, we
can show that 
\[
\left(\forall i\in\{s_{1},\ldots,s_{m}\}\right)\quad f_{s_{2}}(d_{s_{2}}-h_{s_{2}}^{T}z)=f_{s_{i}}(d_{s_{i}}-h_{s_{i}}^{T}z^{*}).
\]
Thus we have arrived at \eqref{eq:goal-pareto-z}.  \hfill  $\blacksquare$

\end{document}